\documentclass[12pt,a4paper]{amsart}

\numberwithin{equation}{subsection}

\newtheorem{theorem}[subsection]{Theorem}
\newtheorem{lemma}[subsection]{Lemma}

\newtheorem{corollary}[subsection]{Corollary}
\newtheorem{conjecture}[subsection]{Conjecture}
\newtheorem{remark}[subsection]{Remark}

\theoremstyle{definition}

\theoremstyle{remark}

\newcommand{\ov}{\overline}
\newcommand{\hilbert}{{\mathcal H}}

\newcommand{\fR}{{\mathfrak R}}

\newcommand{\zp}{{{\bf{Z}}/p}}

\newcommand{\fq}{{\bf{F}}_q}
\newcommand{\fp}{{\bf{F}}_p}
\newcommand{\field}{{\bf{F}}}

\newcommand{\fieldvg}{{\field}[V]^G}
\newcommand{\fieldv}{{\field}[V]}

\newcommand{\ivcoin}{\field[V_4]_{\zp}}
\newcommand{\vcoin}{\field[V_5]_{\zp}}

\newcommand{\tr}{\mathop{\rm Tr}}
\newcommand{\Tr}{\mathop{\rm Tr}}

\newcommand{\lt}{\mathop{\rm LT}}
\newcommand{\lm}{\mathop{\rm LM}}
\newcommand{\wt}{\mathop{\rm wt}}

\newcommand{\del}{\Delta}

\newcommand{\tat}{t\^ete-a-t\^ete}
\newcommand{\binomial}[2]{\genfrac{(}{)}{0pt}{}{#1}{#2}}

\title[Modular Coinvariants]{Coinvariants for modular representations of cyclic groups of prime order}

\author{M\"ufit Sezer}
\author{R.~James Shank}

\address{Institute of Mathematics, Statistics \&  Actuarial Science \\
 \hfil\break\indent University of Kent at Canterbury, CT2 7NF, UK}

\email{M.Sezer@kent.ac.uk}
\email{R.J.Shank@kent.ac.uk}

\thanks{Research supported by a grant from EPSRC}

\date{\today}

\begin{document}

\begin{abstract} We consider the ring of coinvariants for modular representations of cyclic groups of
prime order. For all cases for which explicit generators for the ring of invariants are known,
we give a reduced Gr\"obner basis for the Hilbert ideal and the corresponding monomial basis
for the coinvariants. We also describe the decomposition of the coinvariants as a module over the group ring.
For one family of representations, we are able to describe the coinvariants despite the fact that an explicit
generating set for the invariants is not known. In all cases our results confirm
the conjecture of Harm Derksen and Gregor Kemper on degree bounds for generators of the Hilbert ideal.
As an incidental result, we identify the coefficients of the
monomials appearing in the orbit product of a terminal variable for the three dimensional indecomposable representation. 
\end{abstract}

\maketitle

\section{Introduction}

Let $V$ denote a finite dimensional representation of a finite group $G$ over a field $\field$.
If the characteristic of $\field$ divides the order of $G$, then $V$ in called a {\it modular}
representation. Choose a basis $\{X_1,\ldots,X_n\}$ for the dual vector space $V^*$.  
The action of $G$ on $V$ induces an action on $V^*$ which extends to an
action by algebra automorphisms on the symmetric algebra $\field[V]:=S(V^*)=\field[X_1,\ldots,X_n]$.
The ring of {\it invariants},
$$\field[V]^G:=\{f\in \field[V]\mid g(f)=f,\ \forall g\in G\},$$
is a finitely generated subring of $\field[V]$. The {\it Noether number}, $\beta(V)$, 
is defined to be the least integer $d$ such that
$\fieldv^G$ is generated by homogeneous elements of degree less than or equal to $d$.
 The {\it Hilbert ideal}, which we denote by $\hilbert$, is the ideal in
$\field[V]$ generated by the homogeneous invariants of positive degree and the ring of {\it coinvariants} is the quotient
$$\field[V]_G:=\field[V]/\hilbert.$$
Since the Hilbert ideal is closed under the group action, the coinvariants are a module over the group ring $\field G$.
Furthermore, since $G$ is finite, $\fieldv$ and $\fieldvg$ have the same Krull dimension.
Therefore $\field[V]_G$ is a finite dimensional graded $\field$-algebra.
Let ${\rm td}(\field[V]_G)$ denote the {\it top degree} of $\field[V]_G$, i.e., the largest
degree in which $\field[V]_G$ is non-zero.
The ring of coinvariants has been studied extensively for $\field$ a field of characteristic zero, particularly for $V$ a 
reflection representation. For reflection representations in characteristic zero, the coinvariants are isomorphic, 
as a module over the group ring, to the regular representation
(see, for example, \cite{chevalley}, \cite[Ch.\,V, \S 5.2]{bourbaki} or \cite[Ch.\,VII, \S24-1]{kane}). 
Coinvariants in characteristic zero continue to attract attention (see, for example, \cite{gordon}, \cite{hhlru} and \cite{haiman} ).
Relatively little is known about coinvariants for modular representations. 
The coinvariants for the natural modular representations of ${\rm GL}_n(\fq)$ and its $p$-Sylow subgroup were considered by 
Campbell et al.\ in \cite{chsw}.
Larry Smith has investigated modular coinvariants for two and three dimensional representations \cite{smith:st-proc}
and in the case that the invariants are a polynomial algebra (\cite{smith:st-glag}, \cite{smith:edin}).
In this paper we consider the coinvariants for the simplest modular representations, the modular representations of cyclic groups of prime order.

For the remainder of the paper, let $p$ denote a prime number, let $\zp$ denote the cyclic group of
order $p$ and let $\field$ denote a field of characteristic $p$.
A representation of a cyclic group is determined by the Jordan canonical
form of the image of the generator.  If $n\le p$ then the $n\times n$
matrix over ${\bf F}$ consisting of a single Jordan block with
eigenvalue $1$, has order $p$ and determines an indecomposable
representation of $\zp$ which we denote by $V_n$ (For $n>p$, the order
of the matrix is greater than $p$.).  Note that there are no non-trivial
$p^{\rm th}$ roots of unity in ${\bf F}$.  Thus $1$ is the only
eigenvalue for the image of a generator of $\zp$ under a
representation over ${\bf F}$. Therefore, up to isomorphism, the only indecomposable
${\bf F}\zp$ -- modules are $V_1, V_2,\ldots, V_p$. We will denote the
direct sum of $m$ copies of $V_n$ by $mV_n$.

Despite the simplicity of the representation theory, computing explicit generators for
$\fieldv^\zp$ is a relatively difficult problem. 
 Minimal generating sets for
${\bf F}[V_2]^{\zp}$ and ${\bf F}[V_3]^{\zp}$ can be found in
Dickson's Madison Colloquium \cite{dickson}.
Finite SAGBI bases\footnote{A SAGBI basis is a particularly nice generating set.}
for ${\bf F}[V_4]^{\zp}$ and ${\bf F}[V_5]^{\zp}$ can be found in
\cite{shank}. The problem of finding an explicit generating set for
${\bf F}[V_n]^{\zp}$ for $n>5$ remains open.
Even when the invariants of the indecomposable summands are understood, it can
be difficult to construct generating sets for decomposable
representations. Campbell \& Hughes, in \cite{ch},
describe a generating set for ${\bf F}[mV_2]^{\zp}$ which is refined to
a minimal generating set in \cite{sw:comp}.
SAGBI bases are given for ${\bf F}[V_2\oplus V_3]^{\zp}$ in
\cite{sw:noeth} and $\field[2V_3]^{\zp}$ in \cite{cfw}.
We refer to an $\field\zp$ -- module as {\it reduced} if it is
a direct sum of non-trivial modules. In summary,
the only reduced representations for which explicit generating sets for the ring of invariants are known 
are: $mV_2$, $V_2\oplus V_3$, $V_3$, $2V_3$, $V_4$, $V_5$.
For each of these representations we will give a reduced Gr\"obner basis for the Hilbert ideal
and describe the corresponding monomial basis for the coinvariants. We will also use the 
monomial basis to describe the $\field\zp$ -- module structure of the coinvariants.
By relating $mV_2\oplus\ell V_3$ to $(m+\ell)V_2$, we are able describe $\field[mV_2\oplus\ell V_3]_{\zp}$
despite the fact that an explicit generating set is not known for $\field[mV_2\oplus\ell V_3]^{\zp}$.
Our results give $(m+\ell)(p-1)+1$ as an upper bound on the degrees of a minimal generating
set for $\field[mV_2\oplus\ell V_3]^{\zp}$.  Harm Derksen and Gregor Kemper have conjectured
that the order of the group is an upper bound on the degrees of a minimal homogeneous generating set for
the Hilbert ideal \cite[3.8.6 (b)]{dk}. For all of the examples considered here, our calculations confirm
this conjecture. We note that  $\field[2V_2]_{{\bf Z}/2}$ was considered in \cite{smith:st-proc}.

Let $\sigma$ denote a generator of $\zp$.  In the group ring
$\field\zp$, define $\del:=\sigma-1$ and $\tr:=\sum_{i=1}^p \sigma^i$.
The kernel of $\del$ acting on a module gives the invariant elements
in the module and $\tr$ gives a homomorphism of $\fieldv^{\zp}$ --
modules from $\fieldv$ to $\fieldv^{\zp}$ known as the {\it transfer}.
The image of the transfer is an ideal in the ring of
invariants. Observe that a basis for the coinvariants lifts to a set
of generators for $\fieldv$ as a module over $\fieldv^{\zp}$.
Applying the transfer to a set of module generators gives a generating
set for the image of the transfer as an ideal. Thus a basis for the
coinvariants gives a generating set for the image of the transfer and
the largest degree of a basis element gives an upper bound on the
degrees of a generating set for the image of the transfer. 
It is a consequence of  \cite[4.2 \& 6.3]{sw:noeth} and \cite[4.1]{shank}
that for $n>3$, ${\rm td}(\field[V_n]_{\zp})\geq \beta(V_n)\geq 2p-3$. 
The results in this paper support the following
strengthening of \cite[Conjecture 6.1]{sw:noeth}.

\begin{conjecture} For $n>3$, ${\rm td}(\field[V_n]_{\zp})=2p-3$.
\end{conjecture}

For an element $\varphi\in V^*$, define the {\it norm} of $\varphi$ to be the
product over the orbit of $\varphi$. Thus, if $\varphi\in V^*\setminus (V^*)^{\zp}$, 
$N(\varphi):=\prod_{i=1}^p\sigma^i(\varphi)$. If we choose a basis $\{X,Y,Z\}$ for $V_3^*$
so that $\del(Z)=Y$, $\del(Y)=X$ and $\del(X)=0$, then $\field[V_3]^{\zp}$ is the
hypersurface generated by $X$, $Y^2-X(Y+2Z)$, $N(Y)$ and $N(Z)$.
It is well known that $N(Y)=Y^p-YX^{p-1}$. However, the expansion of $N(Z)$ is far more complicated
and, to our knowledge, does not appear in the literature.
Knowledge of certain coefficients in the expansion was necessary for some of our calculations.
In Section~\ref{normsection}, we have worked out a complete description of the expansion.

We adopt the convention of using upper case letters to denote variables in $\field[V]$ and
the corresponding lower case letters to denote the images of the variables in $\field[V]_G$.
We use the term {\it monomial} to mean a product of variables. For an ideal $I$,
we write $f\equiv_I h$ if $f-h\in I$.
As a general reference for the invariant theory of finite groups see
Benson \cite{benson}, Derksen \& Kemper \cite{dk}, Neusel \& Smith \cite{ns} or Smith \cite{smith:book}.
As a reference for Gr\"obner bases we recommend Adams \& Loustaunau \cite{al} or Sturmfels \cite{sturm}.

\section{The expansion of  $N(Z)$}\label{normsection}\label{expsec}

In this section we describe the expansion of the norm 
of an $\field\zp$ -- module generator of $V_3^*$. 
Choose a basis $\{X,Y,Z\}$ for $V_3^*$ with $\del(Z)=Y$, $\del(Y)=X$ and
$\del(X)=0$. Write $N(Z)=A_0+A_1X+\cdots+A_pX^p$ with each $A_i\in\field[Y,Z]$.

\begin{theorem}\label{expthm} $A_0=Z^p-ZY^{p-1}$, $A_p=A_{p-1}=0$ and
$$A_i=
\begin{cases}
\sum_{k=1}^{i+1}\xi_{ik}Z^kY^{p-i-k} & \text{ for  $1\le i\le \frac{p-1}{2}$,}\cr 
\sum_{k=1}^{p-i}\xi_{ik}Z^kY^{p-i-k} & \text{ for  $\frac{p+1}{2}\le i\le p-2$,}
\end{cases}
$$
where $\xi_{ik}=\frac{(-1)^{i}}{2^i(p-k)}{p-2k+1\choose i-k+1}{p-k\choose k-1}$.
\end{theorem}

The proof of Theorem~\ref{expthm} follows Lemma~\ref{norm3}.
We start with a number of combinatorial lemmas
concerning $\fp$, the field with $p$ elements. The first lemma is well known.

\begin{lemma}
\label{j4}
For a positive integer $\ell$, 
$$\sum_{t\in  \bf F_p}t^{\ell}=
\begin{cases}
-1 &\text { if $p-1$ divides $\ell$;} \cr \\
0 &\text{ if $p-1$ does not divide $\ell$}.
\end{cases}\
$$
\end{lemma}
\begin{proof}
See, for example, \cite[9.4]{chsw}.
\end{proof}

Let $S_i$ denote the set of subsets of ${\bf F}_p$ of size $i$ and, for $j\in {\bf F}_p$,   
let $S_{i,j}$ denote the set of subsets of   ${\bf F}_p$ of size $i$  not containing $j$. 
For  $\alpha \subseteq  {\bf F}_p$, let $\sigma_k(\alpha)$ denote the $k^{\rm th}$ elementary symmetric polynomial 
in the elements of $\alpha$. For convenience, we set $\sigma_0(\alpha)=1$ and to simplify notation we will 
denote $\sigma_i(\alpha)$ by $\pi(\alpha)$ for  $\alpha \in S_i$.
For $j\leq k$, define functions $b_{k,j}:  {\bf F}_p  \rightarrow  {\bf F}_p$ by 
$$b_{k,j}(t):=\sum_{\alpha \in S_{k-1,t}}t\pi(\alpha)\sigma_j(\alpha\cup \{t\})$$
and set $d_{k,j}:=\sum_{\alpha \in S_k}\pi(\alpha)\sigma_j(\alpha)$. Note that $d_{0,0}=1$.

\begin{lemma}
\label{bas}
(i) $\sum_{i\in {\bf F}_p}b_{k,j}(i)=kd_{k,j}$.\\
(ii) $d_{k,j}=b_{k,j}(t)+\sum_{\alpha \in S_{k,t}}\pi(\alpha)\sigma_j(\alpha).$
\end{lemma}
\begin{proof}
The first statement follows from the fact that each term of $d_{k,j}$ appears $k$ times in $\sum_{i\in {\bf F}_p}b_{k,j}(i)$. 
The second statement follows from partitioning $S_k$ into subsets with $t$ and subsets without $t$.
\end{proof}

\begin{lemma}
\label{anina}
For $1\le k < p$, $b_{k,0}(t)=(-1)^{k+1}t^k$ and 
$$d_{k,0}=
\begin{cases}
0 &\text { if } k<p-1; \cr \\
-1 &\text { if } k=p-1.
\end{cases}\
$$ Furthermore $b_{p,0}(t)=t^p-t$ and $d_{p,0}=0$.
\end{lemma}
\begin{proof}
The value of $d_{k,0}$ follows from the description of $b_{k,0}(t)$ using Lemmas \ref{j4} and \ref{bas}.
We prove the given formula for $b_{k,0}(t)$ by induction on $k$.
 Since the product over the empty set is $1$, we have $b_{1,0}(t)=t$.  
Using Lemma \ref{bas}, we see that 
$$ \sum_{\alpha \in S_{k,t}}\pi(\alpha)  =d_{k,0}-t\sum_{\alpha \in S_{k-1,t}}\pi(\alpha)=d_{k,0}-b_{k,0}(t).$$
Thus $b_{k+1,0}(t)=td_{k,0}-tb_{k,0}(t)$. For $k<p-1$, the induction  hypothesis gives $d_{k,0}=0$. 
Therefore $b_{k+1,0}(t)=-t((-1)^{k+1}t^k)=(-1)^{k+2}t^{k+1}$ as required. For $k=p-1$, 
$b_{p,0}(t)=td_{p-1,0}-tb_{p-1,0}t=-t-t((-1)^pt^{p-1})=t^p-t.$ 
\end{proof}

\begin{lemma}\label{sdlem}
For $1\leq k+j<p$, with $0\leq j\le k$, $b_{k,j}(t)=(-1)^{k+1}{k\choose j}t^{k+j}$ and 
$$d_{k,j}=
\begin{cases}
0 &\text { if } k+j<p-1; \cr \\
(-1)^k{k\choose j}\frac{1}{k} &\text { if } k+j=p-1.
\end{cases}\
$$
\end{lemma}
\begin{proof}
The values of $d_{k,j}$ follow from description of $b_{k,j}(t)$ using Lemmas \ref{j4} and  \ref{bas}. 
We prove the formula for $b_{k,j}(t)$ by induction on $k+j$.
For $j=0$ we have $b_{k,0}(t)=(-1)^{k+1}t^k$ by Lemma~\ref{anina}. Working directly from the definition,
$b_{1,1}(t)=t^2$.
Thus the formula holds for $k+j=1$ and $k+j=2$.
Expanding the second factor of each term gives
\begin{eqnarray*}
 b_{k+1,j}(t)&=&\sum_{\alpha \in S_{k,t}}t\pi(\alpha)\sigma_j(\alpha\cup \{t\})\\
             &=&\sum_{\alpha \in S_{k,t}}t\pi(\alpha)t\sigma_{j-1}(\alpha)+\sum_{\alpha \in S_{k,t}}t\pi(\alpha)\sigma_j(\alpha)\\
             &=& t^2\sum_{\alpha \in S_{k,t}}\pi(\alpha)\sigma_{j-1}(\alpha)+t\sum_{\alpha \in S_{k,t}}\pi(\alpha)\sigma_j(\alpha).
\end{eqnarray*}
So by  Lemma~\ref{bas}(ii), we have 
$$b_{k+1,j}(t)=t^2(d_{k,j-1}-b_{k,j-1}(t))+t(d_{k,j}-b_{k,j}(t)).$$
For $2\leq k+j<p-1$, the induction hypothesis gives $d_{k,j-1}=d_{k,j}=0$. 
Therefore  
\begin{eqnarray*}
b_{k+1,j}(t)&=&t^2(-1)^{k+2}{k \choose j-1}t^{k+j-1}+t(-1)^{k+2}{k\choose j}t^{k+j}\\
            &=&(-1)^{k+2}t^{k+j+1}\left({k\choose j-1}+{k\choose j}\right)=(-1)^{k+2}t^{k+j+1}{k+1\choose j},
\end{eqnarray*}
as desired.
\end{proof}

\begin{lemma}\label{ldlem}
Suppose $p-1<k+j<2p-2$. Then $b_{k,j}(t)=(-1)^{k+1}{k\choose j}t^{k+j}+f(t)$, 
where $f(t)$ is a polynomial of degree less or equal to $k+j-(p-1)$, and $d_{k,j}=0$.
\end{lemma}
\begin{proof}
The values of  $d_{k,j}$ follow from the description of $b_{k,j}(t)$ using Lemmas \ref{j4} and  \ref{bas}.
The proof of the formula for $b_{k,j}(t)$ is by induction on $k+j$. 
We use the recursive relation from the proof of the previous lemma,
$$b_{k+1,j}(t)=t^2(d_{k,j-1}-b_{k,j-1}(t))+t(d_{k,j}-b_{k,j}(t)).$$
For $j+k=p$, this gives
\begin{eqnarray*}
b_{k,j}(t)&=&t^2(-1)^{k+1}{k-1\choose j-1}t^{k+j-2}+t(-1)^{k-1}\left({k-1\choose j}\frac{1}{k-1}+{k-1\choose j}t^{k+j-1}\right)\\
          &=&(-1)^{k+1}\left({k-1\choose j-1}t^{k+j}+{k-1\choose j}t^{k+j}+t{k-1\choose j}\frac{1}{k-1}\right)\\
          &=&(-1)^{k+1}{k\choose j}t^{k+j}+t{k-1\choose j}\frac{(-1)^{k-1}}{k-1}.
\end{eqnarray*}
For $j+k=p+1$, the recursive relation gives
\begin{eqnarray*}
b_{k,j}(t)&=&t^2(-1)^{k+1}{k-1 \choose j-1}\left(\frac 1 {k-1}+t^{k+j-2}\right)+t(-1)^{k+1}{k-1 \choose j}t^{k-1+j}\\
          &=&(-1)^{k+1}{k \choose j}t^{j+k}+ t^2{k-1 \choose j-1}\frac{(-1)^{k-1}}{k-1}.                        
\end{eqnarray*}
For $p+1<k+j<2p-2$, the induction hypothesis gives $d_{k,j-1}=d_{k,j}=0$. Therefore 
$b_{k+1,j}(t)= t^2((-1)^{k+2}{k\choose j-1}t^{k+j-1}+p(t))+t((-1)^{k+2}{k\choose j}t^{k+j}+q(t))$, 
where $p(t)$ is a polynomial of degree less than or equal to $k+j-1-(p-1)$ and 
$q(t)$ is a polynomial of degree less than or equal to $k+j-(p-1)$.
Collecting terms gives 
$b_{k+1,j}(t)
=(-1)^{k+2}t^{k+j+1}{k+1\choose j}+t^2p(t)+tq(t)$. 
Since $t^2p(t)+tq(t)$ is a polynomial of degree at most $k+j+1-(p-1)$, the result follows.
\end{proof}

For a set $\gamma \subseteq {\bf F}_p$,  let $S_{b,\gamma}$  denote the set of subsets of ${\bf F}_p$ of size 
$b$ that do not contain any element from the set $\gamma$. We note a counting lemma.

\begin{lemma}\label{count}
$$\sum_{\gamma \in S_c,\ \alpha \in S_{b,\gamma}}\sigma_j(\gamma)\pi(\gamma)\pi(\alpha)={b+c-j\choose b}d_{b+c,j}.$$
\end{lemma}
\begin{proof}
Recall that $d_{b+c,j}=\sum_{\theta\in S_{b+c}}\pi(\theta)\sigma_j(\theta)$. Each term in $\pi(\theta)\sigma_j(\theta)$ 
is of the form $\pi(\tau)\pi(\theta)$ for $\tau$ a subset of $\theta$ of size $j$. 
The term $\pi(\tau)\pi(\theta)$ occurs ${b+c-j\choose b}$ times on the left hand side of the equation, 
once for each choice of $\alpha\in \theta\setminus \tau$. 
\end{proof}

 Let $A_{b,c}$ denote the coefficient   of $X^cY^bZ^{p-c-b}$ in $N(Z)$.

\begin{lemma}
\label{norm3}
(i) Suppose $0<c<p-1$. If there exists an integer $j$ satisfying $0\leq j\leq c$ and $b+c+j=p-1$, then
$$A_{b,c}= \frac{(-1)^{b+2c-j}{b+c-j\choose c-j}{b+c\choose j}}{2^c(b+c)};$$
otherwise $A_{b,c}=0$.\\
(ii) $A_p=A_{p-1}=0$.\\
(iii) $A_0=Z^p-ZY^{p-1}$.
\end{lemma}
\begin{proof}
Recall that $\sigma^m(Z)=Z+mY+{m \choose 2}X$.
By identifying the terms in $\prod_{m\in\fp}\sigma^m(Z)$ which contribute to
the coefficient of $X^cY^bZ^{p-c-b}$ we see that
\begin{eqnarray*}
A_{b,c}&=&\sum_{\{i_1, \dots , i_c\}\in S_c}\ \sum_{\{j_1,\dots j_b\}\in S_{b,\{i_1, \dots , i_c\}}}
           {i_1 \choose 2}{i_2 \choose 2}\cdots {i_c\choose 2}j_1\cdots j_b\\
       &=&\frac{1}{2^c}\sum_{\gamma\in S_c,\ \alpha\in S_{b,\gamma}}\pi(\alpha)\prod_{i\in \gamma}(i^2-i).
\end{eqnarray*}
Expanding gives
\begin{eqnarray*}
\prod_{i\in \gamma}(i^2-i)&=&\sum_{\beta \subseteq \gamma}(-1)^{|\gamma\setminus \beta|}\pi(\beta)\pi(\gamma) \\
                          &=&\pi(\gamma)\sum_{\ell=0}^c (-1)^{c-l}\sigma_{\ell}(\gamma).
\end{eqnarray*}
Substituting this into the previous expression gives
\begin{eqnarray*}
A_{b,c}&=&\frac{1}{2^c}\sum_{\gamma\in S_c,\ \alpha\in S_{b,\gamma}}\pi(\alpha)\pi(\gamma)\sum_{\ell=0}^c(-1)^{c-l}\sigma_{\ell}(\gamma)\\
&=&\frac{1}{2^c}\sum_{\ell=0}^c(-1)^{c-\ell}\left(\sum_{\gamma\in S_c,\ \alpha\in S_{b,\gamma}}\sigma_{\ell}(\gamma)\pi(\alpha)\pi(\gamma)\right).
\end{eqnarray*}
Using Lemma~\ref{count}, gives 
$$A_{b,c}=\frac{1}{2^c}\sum_{\ell=0}^c (-1)^{c-\ell}{b+c-\ell\choose b}d_{b+c,\ell}.$$

It follows from the definition of $d_{k,j}$ that $d_{p,j}=0$. Thus if $b+c=p$, $A_{b,c}=0$. Therefore we may assume $c\leq p-1$.
Using Lemma~\ref{ldlem} and Lemma~\ref{sdlem}, if $0<c+b+j<2p-2$ then $d_{b+c,j}=0$ unless $b+c+j=p-1$. 
If $c=p-1$ and $b+c<p$, then $b=0$. In this case the above summation gives $A_{0,p-1}=\frac{1}{2^{p-1}}(d_{p-1,0}+d_{p-1,p-1})$.
However, explicit calculation gives $d_{p-1,p-1}=1$ and Lemma~\ref{anina} gives $d_{p-1,0}=-1$. Thus $A_{0,p-1}=0$.
For $c=0$, we have $A_{0,0}=1$, $A_{p-1,0}=-1$ and all other $A_{b,0}=0$.
For $0<c<p-1$, we have $0<b+c+\ell<2p-2$. Therefore, there is at most one non-zero term in the above summation. 
If there exists $j\leq c$ with $b+c+j=p-1$ then, using Lemma~\ref{sdlem}, there is non-zero term and
\begin{eqnarray*}
A_{b,c}&=&\frac{(-1)^{c-j}}{2^c} {b+c-j\choose b}d_{b+c,j}\\
       &=&\frac{(-1)^{c-j}}{2^c} {b+c-j\choose b}(-1)^{b+c}{b+c\choose j}\frac{1}{b+c}\\
       &=&\frac{(-1)^{b+2c-j}}{2^c(b+c)}{b+c-j \choose b}{b+c\choose j}
\end{eqnarray*}
as required. If no solution exists, $A_{b,c}=0$.
\end{proof}

To complete the proof of Theorem~\ref{expthm}, we need to identify $A_c$ for $0<c<p-1$.
For $k=p-b-c$ and $A_{b,c}\not=0$, we have $k=p-c-(p-1-c-j)=j+1$. Substituting 
$b=p-k-c$ and $j=k-1$ into the formula for $A_{b,c}$ gives
\begin{eqnarray*}
A_{p-k-c,c}&=&\frac{(-1)^{p-k+c-(k-1)}}{2^c(p-k)}{p-k-(k-1) \choose p-k-c}{p-k\choose k-1}\\
           &=&\frac{(-1)^{p-2k+c-1}}{2^c(p-k)}{p-2k+1 \choose c-k+1} =\xi_{k,c}.
\end{eqnarray*}
For a fixed $c$, the summation is from $k=1$ to $k=c+1$ subject to the condition that $k+c\leq p$.
For $c\leq (p-1)/2$, this condition imposes no restriction. For $c\geq (p+1)/2$, the summation terminates with
$k=p-c$.

\section{ The coinvariants of $mV_2\oplus\ell V_3$}

We start by describing the coinvariants of $mV_2$. Choose a basis $\{X_i,Y_i\mid i=1,\ldots,m\}$ for $(mV_2)^*$ with
$\del(Y_i)=X_i$ and $\del(X_i)=0$. We use the graded reverse lexicographic order with
$X_i<Y_i<X_{i+1}$. For $i=1,\ldots,m$ and $i<j$, define $u_{ij}:=X_jY_i-X_iY_j$. Campbell and Hughes \cite{ch}
have shown that $$\{X_i,N(Y_i),u_{ij}\mid i=1,\ldots,m;\ i<j\}\cup\{\tr(\beta)\mid \beta\ {\rm divides}\ (Y_1\cdots Y_m)^{p-1}\}$$
is a generating set for $\field[mV_2]^{\zp}$. It is well known that $N(Y_i)=Y_i^p-Y_iX_i^{p-1}$. Furthermore, if $\beta$ 
divides $(Y_1\cdots Y_m)^{p-1}$, then $\del(\beta)\in (X_1,\ldots,X_m)\field[mV_2]$. Thus 
$\tr(\beta)=\del^{p-1}(\beta)\in (X_1,\ldots,X_m)\field[mV_2]$. As a consequence, we have the following.

\begin{theorem}\label{mv2thm} A reduced Gr\"obner basis for the Hilbert ideal of $mV_2$ is given by $\{X_i,Y_i^p\mid i= 1,\ldots, m\}$,
the corresponding monomial basis for $\field[mV_2]_{\zp}$ is given by the monomial factors of $(y_1\cdots y_m)^{p-1}$,
and $\field[mV_2]_{\zp}$ is a trivial $\field\zp$ -- module.
\end{theorem}

For the rest of this section, we assume $p>2$.
The natural inclusion of $(m+\ell)V_2$ into $mV_2\oplus \ell V_3$ induces an algebra epimorphism 
$\rho:\field[mV_2\oplus \ell V_3]\to \field[(m+\ell)V_2]$. We will use this map in conjunction with Theorem~\ref{mv2thm}
to describe the coinvariants of $mV_2\oplus \ell V_3$. Choose a basis
$$\{X_i,Y_i\mid i=1,\ldots,m\}\cup\{X_i,Y_i,Z_i\mid i=m+1,\ldots,m+\ell\}$$
for $(mV_2\oplus \ell V_3)^*$ with $\del(Z_i)=Y_i$, $\del(Y_i)=X_i$ and $\del(X_i)=0$.
We use the graded reverse lexicographic order with
$X_i<Y_i<Z_i<X_{i+1}$.
For $i=1,\ldots,m+\ell$ and $i<j$, define $u_{ij}:=X_jY_i-X_iY_j$ and, for $i=m+1,\ldots,m+\ell$ and $i<j$,
define $d_i:=Y_i^2-X_i(Y_i+2Z_i)$ and $w_{ij}:=Z_iX_j-Y_iY_j+X_iZ_j+X_iY_j$.
A straightforward calculation verifies that $u_{ij}$, $d_i$ and $w_{ij}$ are all elements of
$\field[mV_2\oplus\ell V_3]^{\zp}$. Let $I$ be the ideal in  $\field[mV_2\oplus\ell V_3]$ generated by
$$\{X_i,N(Y_i)\mid i=1,\ldots m\}\cup\{X_i,d_i,w_{ij},N(Z_i)\mid i=m+1,\ldots,m+l;\ i<j\}.$$
and define
$$\Lambda:=\{X_i, Y_i^p\mid i=1,\ldots,m\}\cup\{X_i, Y_iY_j,Z_i^p\mid i=m+1,\ldots,m+\ell;\ i\leq j\}.$$

\begin{lemma} The set $\Lambda$ is a reduced Gr\"obner basis for $I$.
\end{lemma}
\begin{proof}
It follows from Section~\ref{normsection} that
$N(Z_i)\equiv_{(X_i)}Z_i^p-Z_iY_i^{p-1}$. Using this, along with the
expansion of $N(Y_i)$ given above and the definition of $d_i$ and
$w_{ij}$, it is clear that $\Lambda$ generates $I$. Since $\Lambda$ is
a set of monomials and a minimal generating set for $I$, it is a
reduced Gr\"obner basis for $I$.
\end{proof}

\begin{lemma}\label{trlem} If $\beta$ divides $(Y_1\cdots Y_mZ_{m+1}\cdots Z_{m+\ell})^{p-1}$,
then $\tr(\beta)\in I$.
\end{lemma}

\begin{proof}
Write $\beta=Y^FZ^E$ where $Y^F:=\prod_{i=1}^mY_i^{f_i}$ with $F:=(f_1,\ldots,f_m)\in {\bf Z}^m$
and $Z^E:=\prod_{i=m+1}^{m+\ell}Z_i^{e_i}$ with $E:=(e_{m+1},\ldots,e_{m+\ell})\in {\bf Z}^{\ell}$. 
Clearly $\del(Y_i)\equiv_I 0$. Therefore $\del(\beta)=Y^F\del(Z^E)$ and
$\tr(\beta)=\del^{p-1}(\beta)=Y^F\tr(Z^E)$. Thus it is sufficient to show that $\tr(Z^E)\in I$.
Recall that $\sigma^c(Z_i)=Z_i+cY_i+\binomial c 2 X_i\equiv_I Z_i+cY_i$. Thus
\begin{eqnarray*}
\tr(Z^E) &=&\sum_{c\in\fp }\sigma^c(Z^E)\\
         &\equiv_I& \sum_{c\in\fp}\prod_{i=m+1}^{m+\ell}\left(Z_i+cY_i\right)^{e_i}.
\end{eqnarray*}
Using the fact that, for $i=m+1,\ldots,m+\ell$, we have $Y_i^2\in I$, gives
$$\tr(Z^E)\equiv_I \sum_{c\in\fp}\prod_{i=m+1}^{m+\ell}\left(Z_i^{e_i}+e_icY_iZ_i^{e_i-1}\right).$$
Furthermore, for $i=m+1,\ldots,m+\ell$ and $i<j$, we have $Y_iY_j\in I$. Thus
\begin{eqnarray*}
\tr(Z^E)&\equiv_I& \sum_{c\in\fp} \left(Z^E+c\sum_{i=m+1}^{m+\ell}e_iY_i\frac{Z^E}{Z_i}\right)\\
        &\equiv_I& Z^E\left(\sum_{c\in\fp}1\right)+\left(\sum_{c\in\fp}c\right)\left(\sum_{i=m+1}^{m+\ell}e_iY_i\frac{Z^E}{Z_i}\right).
\end{eqnarray*}
Therefore, using Lemma~\ref{j4}, $\tr(Z^E)\equiv_I 0$, as required.
\end{proof}

The algebra epimorphism $\rho:\field[mV_2\oplus\ell V_3]\to\field[(m+\ell)V_2]$, introduced above, 
is a morphism of $\field\zp$ -- modules and is determined by 
$\rho(Z_i)=Y_i$, $\rho(Y_i)=X_i$ and $\rho(X_i)=0$ for $i>m$ and by
$\rho(Y_i)=Y_i$ and $\rho(X_i)=X_i$ for $i\leq m$. The kernel of $\rho$ is generated by $\{X_i\mid i=m+1,\ldots,m+\ell\}$ 
and is contained in $I$. Since $\rho$ is surjective, the image of $I$ under $\rho$ is an ideal in $\field[(m+\ell)V_2]$. 
Intersecting this ideal with
the ring of invariants gives the ideal 
$J:=\rho(I)\cap \field[(m+\ell)V_2]^{\zp}$.

\begin{lemma}\label{projlem} The natural projection  from $\field[(m+\ell)V_2]^{\zp}$ to
$\field[(m+\ell)V_2]^{\zp}/J$ induces an epimorphism of vector spaces from
$${\rm Span}\left( \{X_i\mid i=m+1,\ldots,m+\ell\}\cup\{u_{ij}\mid i=1,\ldots,m+\ell;\ i<j\ {\rm and}\ m<j\}\right)$$
to $\field[(m+\ell)V_2]^{\zp}/J$. 
\end{lemma}
\begin{proof}
Recall that $\field[(m+\ell)V_2]^{\zp}$ is generated by
 $$\{X_i,N(Y_i),u_{ij}\mid i=1,\ldots,m+\ell;\ i<j\}\cup\{\tr(\alpha)\mid \alpha\ {\rm divides}\ (Y_1\cdots Y_{m+\ell})^{p-1}\}.$$
For each monomial $\alpha$ dividing $(Y_1\cdots Y_{m+\ell})^{p-1}$, there exists a monomial 
$\beta$ dividing $(Y_1\cdots Y_mZ_{m+1}\cdots Z_{m+\ell})^{p-1}$
with $\rho(\beta)=\alpha$. By Lemma~\ref{trlem}, $\tr(\beta)\in I$. Therefore $\tr(\alpha)=\rho(\tr(\beta))\in J$.
For $i\leq m$, $\rho(Y_i)=Y_i$. Thus $N(Y_i)=\rho\left(N\left(Y_i\right)\right)\in J$. For $i>m$, $\rho(Z_i)=Y_i$ giving
$N(Y_i)=\rho(N(Z_i))\in J$. For $i\leq m$, $X_i=\rho(X_i)\in J$. For $i>m$, $X_i^2=\rho(d_i)\in J$ and $X_iX_j=-\rho(w_{ij})\in J$.
For $i<j\leq m$, $u_{ij}=\rho(u_{ij})\in J$. We have shown that, for all $i$ and $j$, $X_i^2$ and $X_iX_j$ lie in $\rho(I)$.
Therefore $u_{ij}u_{rs}=X_jX_sY_iY_r-X_iX_sY_jY_r-X_jX_rY_iY_s+X_iX_rY_jY_s$ and $X_iu_{rs}=X_iX_sY_r-X_iX_rY_s$ lie in
$\rho(I)$. Since these elements are invariant, they lie in $J$.
\end{proof}

\begin{theorem}\label{mv2lv3thm}
The ideal $I$ coincides with the Hilbert ideal of $mV_2\oplus \ell V_3$.
\end{theorem}
\begin{proof} By definition, $I\subseteq {\mathcal H}$. Thus it is sufficient
to show that every invariant lies in $I$. Suppose that $f$ is a homogeneous element of
$\field[mV_2\oplus \ell V_3]^{\zp}$ with $\deg(f)>2$. Then using Lemma~\ref{projlem},
$\rho(f)\in J\subseteq \rho(I)$.
Thus there exist $\widetilde f\in I$ with $\rho(\widetilde f)=\rho(f)$. Therefore
$\widetilde f - f \in\ker(\rho)\subseteq I$. Thus $f\in I$ as required.

Every homogeneous invariant of degree $1$ is a linear combination of the $X_i$ and hence
lies in $I$. Therefore we need only verify that all homogeneous invariants of degree $2$
lie in $I$. To do this we grade $\field[mV_2\oplus \ell V_3]$ over
${\bf Z}^{m+\ell}=\oplus_{i=1}^{m+\ell}b_i{\bf Z}$ by defining
the multidegree of $X_i$, $Y_i$ and $Z_i$ to be $b_i$. The group action preserves multidegree.
Therefore we may restrict to invariants which are homogeneous with respect to multidegree.
Since the total degree is $2$, the possible multidegrees are $2b_i$ and $b_i+b_j$.
For multidegree $2b_i$, we use the descriptions of
$\field[V_2]^{\zp}$  and $\field[V_3]^{\zp}$ from \cite{dickson}.
For multidegree $b_i+b_j$, we use the description of
$\field[2V_2]^{\zp}$ from \cite{ch}, the description of
$\field[2V_3]^{\zp}$ from \cite{cfw} and the description
of $\field[V_2\oplus V_3]^{\zp}$ from \cite{sw:noeth}. In all cases, the only generators in
degrees less than or equal to $2$ are $X_i$, $d_i$, $u_{ij}$ and $w_{ij}$.
All of these invariants appear in $I$.
\end{proof}

\begin{corollary}\label{mv2lv3cor}
 A reduced Gr\"obner basis for the Hilbert ideal of $mV_2\oplus\ell V_3$ is given by 
$$\{X_i, Y_i^p\mid i=1,\ldots,m\}\cup\{X_i, Y_iY_j,Z_i^p\mid i=m+1,\ldots,m+\ell;\ i\leq j\},$$
the corresponding monomial basis for $\field[mV_2\oplus\ell V_3]_{\zp}$ is given by 
the monomial factors of 
$y_j(y_1\cdots y_m z_{m+1}\cdots z_{m+\ell})^{p-1}$ for $j=m+1,\ldots,m+\ell$, and the Hilbert series of
 $\field[mV_2\oplus\ell V_3]_{\zp}$ is $\left(\ell t+1\right)\left(1+t+\cdots+t^{p-1}\right)^{m+\ell}$.
Furthermore, both as $\field$-algebras and $\field\zp$ -- modules,
$\field[mV_2\oplus\ell V_3]_{\zp}\cong\field[mV_2]_{\zp}\otimes\field[\ell V_3]_{\zp}$.
\end{corollary}

\begin{remark} We have shown that the Hilbert ideal of $mV_2\oplus\ell V_3$ is generated 
by homogeneous elements of degree less than or equal to $p$, the order of the group, confirming
the conjecture of Derksen \& Kemper \cite[3.8.6(b)]{dk} in this case.
Theorem~\ref{h4} and Theorem~\ref{v5thm} confirm the conjecture for
$V_4$ and $V_5$ respectively.
\end{remark}

\begin{corollary} If $m+\ell>2$, then $$(m+\ell)(p-1)\leq\beta(mV_2\oplus\ell V_3)\leq(m+\ell)(p-1)+1.$$
\end{corollary}
\begin{proof}
From \cite[4.2]{sw:noeth}, we know that the Noether number of a representation is greater
than or equal to the Noether number of a subrepresentation. Thus
$\beta\left(\left(m+\ell\right)V_2\right)\leq\beta\left(mV_2\oplus\ell V_3\right)$.
From \cite{ch} or \cite{richman}, for $m+\ell>2$, the Noether number of $(m+\ell)V_2$ is $(m+\ell)(p-1)$.
This gives the first inequality. The second inequality follows from \cite[2.12]{hk} using the fact
that ${\rm td}\left(\field[mV_2\oplus\ell V_3]_{\zp}\right)=(m+\ell)(p-1)+1$ is an upper
bound on the degrees of the generators of the image of the transfer.
\end{proof}

\begin{remark}
The generating sets for $\field[V_2\oplus V_3]^{\zp}$ and $\field[2V_3]^{\zp}$ in \cite{sw:noeth} and \cite{cfw}
respectively, include elements of degree $2(p-1)+1$. However, these generating sets are not proven to be minimal.
MAGMA \cite{magma} calculations for the primes $3$, $5$ and $7$ do give $2(p-1)+1$ as the Noether number for
these representations. Further MAGMA calculations show that $2V_2\oplus V_3$, $V_2\oplus 2V_3$ and $3V_3$ at $p=3$,
all have Noether number $7$. 
\end{remark}

In order to describe the $\field\zp$ -- module structure of $\field[mV_2\oplus\ell V_3]_{\zp}$,
we use the grading introduced in the proof of Theorem~\ref{mv2lv3thm}. Since ${\mathcal H}$ is
generated by elements which are homogeneous with respect to multidegree, the grading
on $\field[mV_2\oplus\ell V_3]$ induces a grading on  $\field[mV_2\oplus\ell V_3]_{\zp}$.
The group action preserves the multidegree. Therefore the homogeneous components give
an $\field\zp$ -- module decomposition.  Furthermore, since 
$\field[mV_2\oplus\ell V_3]_{\zp}\cong\field[mV_2]_{\zp}\otimes\field[\ell V_3]_{\zp}$ and
$\field[mV_2]_{\zp}$ is a trivial $\field\zp$ -- module, it is sufficient
to describe the module structure of $\field[\ell V_3]_{\zp}$.
Using the notation from the proof of Lemma~\ref{trlem}
we can describe the basis elements for $\field[\ell V_3]_{\zp}$  as $y_j^{\varepsilon}z^E$ 
where $j>m$, $\varepsilon\in\{0,1\}$ and $E=(e_1,\ldots,e_{\ell})\in {\bf Z}^{\ell}$ with
$0\leq e_i\leq p-1$. It is clear that $\del(y_jz^E)=0$ and
$$\del(z^E)=\sum_{j\in \{m+i\mid e_i\not=0\}}y_j\frac{z^E}{z_j}.$$
Sorting the basis elements into their multidegree components gives the following.

\begin{theorem}
In top degree, $\ell(p-1)+1$, the  $\ell$ multidegree components of $\field[\ell V_3]_{\zp}$
are one dimensional with each component given
by $y_jz^{(p-1,p-1,\ldots,p-1)}\field$. For total degree greater than zero and less than $\ell(p-1)+1$,
each multidegree component is given by the span of $\{z^E,y_jz^E/z_j\mid e_{j-m}\not =0\}$
and is isomorphic to $V_2\oplus (k-1)V_1$ where $k$ is the number of non-zero entries in $E$.
\end{theorem}

\section{The coinvariants of $V_4$}\label{v4coin}

In this section we use the generating set for $\field[V_4]^{\zp}$ given in \cite{shank} to
construct a reduced Gr\"obner basis for the Hilbert ideal.
Choose a basis $\{X_1,X_2,X_3,X_4\}$ for $V_4^*$ with $\del(X_i)=X_{i-1}$ for $i>1$ and $\del(X_1)=0$.
We use the graded reverse lexicographic  order with $X_1<X_2<X_3<X_4$. 
We start with a useful lemma.

\begin{lemma}
\label{12}
Suppose $\beta=X_2^iX_3^j$. Further suppose that $\alpha$ is a monomial with $\alpha<\beta$ and $\deg(\alpha)=\deg(\beta)$.
Then $\alpha$ lies in the ideal generated by $\{X_1, X_2^{i+1}\}$.
\end{lemma}
\begin{proof}
When comparing $\alpha$ and $\beta$ using the graded reverse lexicographic order, we first compare the
exponents of $X_1$ and then, if necessary, the exponents of $X_2$.
\end{proof}

\begin{theorem}
\label{h4}
A reduced Gr\"obner basis for the Hilbert ideal of $V_4$ is given by
$\{X_1,X_2^2,X_2X_3^{p-3},X_3^{p-1},X_4^p\}$, the corresponding monomial basis for $\field[V_4]_{\zp}$
is given by the monomial factors of $x_3^{p-2}x_4^{p-1}$ and $x_2x_3^{p-4}x_4^{p-1}$, and the Hilbert
series of $\field[V_4]_{\zp}$ is given by $(1 +2(t+t^2+\cdots +t^{p-3})+t^{p-2})(1+t+\cdots +t^{p-1})$.
\end{theorem}
\begin{proof}
By  \cite[4.1]{shank},  the ring of invariants is generated by $X_1$, $X_2^2-X_1(X_2+2X_3)$, 
$X_2^3+X_1^2(3X_4-X_2)-3X_1X_2X_3$, $g=X_2^2X_3^2+\cdots$, $N(X_4)$ and the following families:\\
(i) $\Tr (X_3^iX_4^{p-1})$ for $0\le i \le p-2$,\\
(ii) $\Tr (X_3^iX_4^{p-2})$ for $3\le i \le p-2$,\\
(iii) $\Tr (X_4^j)$ for $q\le j \le p-2$,\\
(iv) $\Tr (X_3^2X_4^j)$ for $2l-1\le j \le p-2$.

\noindent
where $l=\frac{p-1}{3}$, $q=2l+1$ if $p\equiv 1$ modulo 3 and $l=\frac{p+1}{3}$,  $q=2l-1$ if $p\equiv -1$ modulo 3.
In the following, we will determine the contribution of each generator to the reduced Gr\"obner basis.
We first note that the ideal generated by  $X_1$, $X_2^2-X_1(X_2+2X_3)$, $X_2^3+X_1^2(3X_4-X_2)-3X_1X_2X_3$  has reduced 
Gr\" obner basis  $\{X_1, X_2^2\}$. Furthermore, by Lemma~\ref{12}, all of the monomials appearing in $g$ 
lie in the ideal $(X_1,X_2^2)$.
 
The leading monomials of the elements in the transfer families above were computed in \cite{shank}. 
Using these results, we compute the contributions  to the reduced the Gr\" obner basis
of the second, third and fourth families. 

For the third family, using \cite[3.2]{shank}, the leading monomials are 
$\lm (\Tr (X_4^j))=X_2^{p-1-j}X_3^{2j-p+1}$ for $q\le j \le p-2$. For $j<p-2$, the leading monomial is divisible by $X_2^2$.
For $j=p-2$, the leading monomial is $X_2X_3^{p-3}$. Using Lemma~\ref{12} all ``non-leading'' monomials are in the ideal
$(X_1,X_2^2)$. Therefore the third family contributes $X_2X_3^{p-3}$ to the reduced Gr\"obner basis.

For the second family of transfers, by  \cite[3.4]{shank} we have  
$\lm (\Tr (X_3^iX_4^{p-2}))=X_2X_3^{i+p-3}$ for $3\le i \le p-2$. Thus each leading monomial is divisible by
$X_2X_3^{p-3}$ and, using Lemma~\ref{12}, the non-leading monomials lie in $(X_1,X_2^2)$. Thus the second family
does not contribute to the reduced Gr\"obner basis. 

For the fourth family,
 by  \cite[3.5]{shank}, we have $\lm (\Tr (X_3^2X_4^j))=X_2^{p-1-j}X_3^{2j-p+3}$ for $2l-1\le j \le p-2$.
For $j<p-2$, the leading monomial is divisible by $X_2^2$. For $j=p-2$ the leading monomial is divisible by
$X_2X_3^{p-3}$. Again using Lemma~\ref{12}, all of the non-leading monomials lie in $(X_1,X_2^2)$.
Therefore the fourth family does not contribute to the reduced Gr\"obner basis.

For the first family, by \cite[3.3]{shank} and \cite[3.2]{shank}, we have
 $\lm (\Tr (X_3^iX_4^{p-1})=X_3^{i+p-1}$  for $0\le i \le p-2$.
Thus the leading monomials are all divisible by $X_3^{p-1}$.
We claim that the non-leading monomials appearing in $\tr(X_3^iX_4^{p-1})$ 
all lie in $(X_1,X_2^2,X_2X_3^{p-3})$.
Therefore, proving the claim will show that the first family contributes $X_3^{p-1}$
to the reduced Gr\"obner basis. To prove the claim, we first observe that
$$\sigma^j(X_3^iX_4^{p-1})= \left(X_3+ jX_2+ {j\choose 2}X_1\right)^i\left(X_4+jX_3+{j\choose 2}X_2+{j\choose 3}X_1\right)^{p-1}.$$
Using Lemma~\ref{j4}, the only term not divisible by $X_1$ or $X_2$ which ``survives'' the summation is $j^{p-1}X_3^{p+i-1}$.
Clearly terms divisible by $X_1$ or $X_2^2$ lie in the ideal $(X_1,X_2^2,X_2X_3^{p-3})$.
Thus we may restrict our attention to monomials of the form $X_2X_3^{p-2+i-a}X_4^a$.
If $p-2+i-a\geq p-3$, this monomial lies in $(X_1,X_2^2,X_2X_3^{p-3})$.
Therefore, it is sufficient to show that if $a>i+1$, the term with monomial $X_2X_3^{p-2+i-a}X_4^a$ 
does not survive the summation.
The coefficient  of $X_2X_3^{p-2+i-a}X_4^a$ in $\sigma^j(X_4^{p-1})$ is 
$(p-1)j^{p-2-a}{ j \choose 2}{p-2 \choose a}+ ij^{p-a}{p-1 \choose a}$.
This coefficient has degree $p-a$ as a polynomial in $j$.
Since $i+1< a$, we  have $p-a<p-(i+1)=(p-1)-i$. Therefore $p-a<p-1$ and, by Lemma~\ref{j4},
the term does not survive the summation, proving the claim.

The only remaining invariant is $N(X_4)$. 
Working modulo $(X_1)$, the variable $X_4$ generates an $\field\zp$ -- module isomorphic to $V_3$.
Thus we may use the results of Section~\ref{expsec}.
Write $N(X_4)\equiv_{(X_1,X_2^2)} A_0+A_1X_2$ for $A_0,A_1\in\field[X_3,X_4]$.
By Theorem~\ref{expthm}, we may take $A_0=X_4^p-X_4X_3^{p-1}$ and
$A_1=\xi_{11}X_4X_3^{p-2} + \xi_{12}X_4^2X_3^{p-3}$. Thus $X_2A_1\in(X_2X_3^{p-3})$
and $N(X_4)-X_4^p\in (X_1,X_2^2,X_2X_3^{p-3},X_3^{p-1})$.
Therefore $N(X_4)$ contributes $X_4^p$ to the reduced Gr\"obner basis.

We have shown that $\{X_1,X_2^2,X_2X_3^{p-3},X_3^{p-1}, X_4^p\}$ generates the
Hilbert ideal. Furthermore, it is clear that this is a minimal generating set
of monomials and is, therefore, a reduced Gr\"obner basis.
The corresponding monomial basis consists of all monomials not divisible by any of the
generators and the description of the Hilbert series comes from the monomial basis.
\end{proof}

\smallskip
\begin{remark}\label{v4rem}
We observe that the top degree of $\field[V_4]_{\zp}$ is $2p-3$.
It is clear that $2p-3$ is an upper bound for the Noether number of
$V_4$.  Using the theory of SAGBI bases it is possible to prove that
that $\tr(X_3^{p-2}X_4^{p-1})$ is indecomposable and, therefore,
$\beta(V_4)=2p-3$. We give a sketch of the proof. For the required
background see \cite{rs} or \cite[Ch.\ 11]{sturm}.

Let ${\mathcal C}$ denote the generating set given above
and define ${\mathcal D}={\mathcal C}\setminus \{\tr(X_3^{p-2}X_4^{p-1})\}$.
Note that the elements of ${\mathcal D}$ all have degree less than $2p-3$.
Recall that ${\mathcal C}$ is a SAGBI basis for $\field[V_4]^{\zp}$.
Therefore ${\mathcal D}$ is ``SAGBI to degree $2p-4$''.
The leading monomial of $\tr(X_3^{p-2}X_4^{p-1})$ is $X_3^{2p-3}$.
The powers of $X_3$ appearing
in $\lm({\mathcal D})$ are $X_3^{p-1},X_3^p,\ldots,X_3^{2p-4}$.
Therefore the leading monomial of $\tr(X_3^{p-2}X_4^{p-1})$ 
does not factor over $\lm({\mathcal D})$ and 
${\mathcal D}$ is not a SAGBI basis for $\field[V_4]^{\zp}$.
Thus either $\tr(X_3^{p-2}X_4^{p-1})$ is indecomposable or
a non-trivial \tat\ from ${\mathcal D}$ subducts to an invariant
with leading monomial $X_3^{2p-3}$. However,
the only monomials in degree $2p-3$ which are greater than
$X_3^{2p-3}$ are  of the form $X_3^{2p-3-a}X_4^a$ and
the only element of ${\mathcal D}$ whose lead monomial is divisible by
$X_4$ is $N(X_4)$.
Therefore the only 
\tat s from ${\mathcal D}$ which could subduct to an invariant with
leading monomial $X_3^{2p-3}$ are of the form $f_1N(X_4)-f_2N(X_4)$. 
However, ${\mathcal D}$ is ``SAGBI to degree $2p-4$''. Therefore the
\tat\ $f_1-f_2$ subducts to zero. Thus  $f_1N(X_4)-f_2N(X_4)$
subducts to zero. Since no \tat\ from ${\mathcal D}$ can subduct to
an invariant with leading monomial $X_3^{2p-3}$, $\tr(X_3^{p-2}X_4^{p-1})$
is indecomposable.
\end{remark}

\section{The coinvariants of $V_5$}\label{v5coin}

In this section we use the generating set for $\field[V_5]^{\zp}$ given in \cite{shank} to
construct a reduced Gr\"obner basis for the Hilbert ideal.
Choose a basis $\{X_1,X_2,X_3,X_4,X_5\}$ for $V_5^*$ with $\del(X_i)=X_{i-1}$ for $i>1$ and $\del(X_1)=0$.
We use the graded reverse lexicographic  order with $X_1<X_2<X_3<X_4<X_5$. 

\begin{theorem}\label{v5thm}
For $p>5$, a reduced Gr\"obner basis for the Hilbert ideal of $V_5$ is given by
$$\{X_1,X_2^2, X_3^2-2X_4X_2-X_3X_2, X_4X_3X_2,X_4^{p-4}X_2, X_4^{p-3}X_3, X_4^{p-1}, X_5^p\},$$        
the corresponding monomial basis for $\field[V_5]_{\zp}$
is given by the monomial factors of $x_4^{p-2}x_5^{p-1}$, $x_3x_4^{p-4}x_5^{p-1}$, $x_2x_4^{p-5}x_5^{p-1}$,
and $x_2x_3x_5^{p-1}$, and the Hilbert
series of $\field[V_4]_{\zp}$ is given by $(1 + 3t+4t^2+3(t^3+ \cdots +t^{p-4})+2t^{p-3}+t^{p-2})(1+t+\cdots +t^{p-1})$.
\end{theorem}

\begin{remark}
For $p=5$, a MAGMA \cite{magma} calculation shows that
 a reduced Gr\"obner basis for the Hilbert ideal of $V_5$ is given by
$$\{X_1,X_2^2,X_3^2-2X_4X_2-X_3X_2, X_2X_3X_4, X_4^2X_3+2X_4^2X_2,X_4^3X_2, X_4^4,X_5^5\},$$ 
the corresponding monomial basis for $\field[V_5]_{{\bf Z}/5}$ is
given by the monomial factors of $x_4^3x_5^4$, $x_3x_4x_5^4$, $x_2x_4^2x_5^4$, and $x_2x_3x_5^4$ and the Hilbert
series of $\field[V_5]_{\zp}$ is given by $(1 +3t+4t^2+2t^3)(1+t+t^2+t^3+t^4)$.
\end{remark}

The rest of the section is devoted to the proof of Theorem~\ref{v5thm}.
The generating set given in \cite[5.1]{shank} consists of a list of
prime independent {\it rational} invariants, a list of
transfers, and $N(X_5)$. The first four rational invariants are
$X_1$, $X_2^2-X_1(X_2+2X_3)$, $X_3^2-X_2(X_3+2X_4)+X_1(X_3+3X_4+2X_5)$ and $X_2^3+X_1^2(3X_4-X_2)-3X_1X_2X_3$.
These invariants contribute $X_1$, $X_2^2$ and $X_3^2-X_2(X_3+2X_4)$ to the reduced Gr\"obner basis.
The fifth rational invariant, denoted  by $\overline{\rm inv}(X_3^3)$ in \cite{shank},
can be computed using the algorithm given in the proof of \cite[2.3]{shank}. Working modulo the ideal generated by
$X_1$, this computation gives $\overline{\rm inv}(X_3^3)\equiv_{(X_1)} 2X_3^3-6X_2X_3X_4+6X_2^2X_5-2X_2^2X_3-3X_2X_3^2-6X_2^2X_4$.
This invariant contributes $X_2X_3X_4$ to the reduced Gr\"obner basis. The sixth rational invariant is in fact decomposable
and was required in \cite[5.1]{shank} in order for the generating set to be a SAGBI basis. Therefore, denoting the ideal generated by
the rational invariants by $\fR$, we have $$\fR=(X_1,X_2^2,X_3^2-X_2(X_3+2X_4),X_2X_3X_4)\fieldv.$$
Note that $X_2X_3^2$ and $X_3^3$ are both elements of $\fR$.

The following lemma will be used in determining the contribution of the image of the transfer to the Hilbert ideal.

\begin{lemma}
\label{tra} Suppose $a$, $b$, $c$ and $d$ are non-negative integers.\\
(i) If $c+2b+3a<p-1$, then  $X_2^aX_3^bX_4^cX_5^d$ does not appear in $\Tr (X_5^i)$.\\
(ii) If\, $i-d+b+2a<p-1$, then  $X_2^aX_3^bX_4^cX_5^d$ does not appear in $\Tr (X_4^kX_5^i)$.
\end{lemma}
\begin{proof}
 Note that $\sigma^j(X_5^i)= (X_5+ jX_4+ {j\choose 2}X_3+{j\choose 3}X_2+{j \choose 4}X_1)^i$. 
Thus the coefficient of  $X_2^aX_3^bX_4^cX_5^d$ in $\sigma^j(X_5^i)$ is 
${i\choose a} {i-a \choose b} {i-a-b \choose c}{j\choose 3}^ a  {j\choose 2}^b j^c$ which is a polynomial of degree  $c+2b+3a$ in $j$. 
Hence by Lemma~\ref{j4} the coefficients will sum to zero under the transfer if  $c+2b+3a<p-1$.

For the second statement, note that   
$$\sigma^j(X_4^kX_5^i)\equiv_{(X_1)}\left(X_4+jX_3+{j\choose 2}X_2\right)^k\left(X_5+jX_4+{j\choose 2}X_3+{j\choose 3}X_2\right)^i.$$ 
We show that the coefficient of $X_2^aX_3^bX_4^cX_5^d$ as a polynomial in $j$ is of degree $2a+b+i-d$. 
Assume that $X_2^{a_1}X_3^{b_1}X_4^{c_1}$ comes from the first factor and  $X_2^{a_2}X_3^{b_2}X_4^{c_2}X_5^d$ comes from the second factor. 
Note that we have $a_1+a_2=a$, $b_1+b_2=b$, $c_1+c_2=c$, $a_1+b_1+c_1=k$, $a_2+b_2+c_2+d=i$. 
The coefficient of $X_2^{a_1}X_3^{b_1}X_4^{c_1}$ in  $\sigma^j(X_4^k)$ is of degree $b_1+2a_1$ in $j$. 
On the other hand the coefficient of $X_2^{a_2}X_3^{b_2}X_4^{c_2}X_5^d$ in  $\sigma^j(X_5^i)$ is of degree $c_2+2b_2+3a_2$ in $j$. 
It follows that the coefficient of the product  $X_2^aX_3^bX_4^cX_5^d$ has degree $c_2+2b_2+b_1+3a_2+2a_1=c_2+b_2+a_2+b+2a=i-d+b+2a$. 
By Lemma \ref{j4} the coefficient will sum to zero under the transfer if $i-d+b+2a<p-1$.
\end{proof}

The generating set in \cite[5.1]{shank} includes one exceptional transfer, $\Tr (X_2X_3X_5^{(p-1)/2})$,
and the following five families:\\
(i) $\Tr (X_4^iX_5^{p-1})$ and $\Tr (X_2X_4^iX_5^{p-1})$ for  $0\le i \le p-2$,\\
(ii) $\Tr (X_4^iX_5^{p-2})$ and $\Tr (X_2X_4^iX_5^{p-2})$ for $3\le i \le p-2$,\\
(iii) $\Tr (X_4^2X_5^i)$ and $\Tr (X_2X_4^2X_5^i)$ for $(p-1)/2\le i \le p-2$,\\
(iv) $\Tr (X_5^i)$ for $(p+1)/2\le i \le p-1$,\\
(v) $\Tr (X_2X_5^i)$ for $(p-1)/2\le i \le p-2$.

We start with the fourth family. By \cite[3.2]{shank} the leading monomial of $\Tr (X_5^i)$ is $X_3^{p-1-i}X_4^{2i-p+1}$. 
Therefore, as $i$ runs from $(p+1)/2$ to $p-1$, the leading monomials are
$X_3^{(p-3)/2}X_4^2,X_3^{(p-5)/2}X_4^4,\ldots,X_3X_4^{p-3},X_4^{p-1}$. The hypothesis $p>5$ means that $(p-3)/2\geq 2$.

First assume $i\leq p-4$. In this case, the leading monomial is divisible by $X_3^3$ and 
hence lies in $\fR$. Suppose $\alpha$ is a monomial of degree $i$ with 
$\alpha<X_3^{p-1-i}X_4^{2i-p+1}$ and $\alpha\not\in\fR$.
Since we are using the graded reverse lexicographic order and $p-1-i\geq 3$,
$\alpha$ must be divisible by $X_1$, $X_2$ or $X_3^3$.
Note that $X_1$, $X_2^2$, $X_2X_3^2$, $X_3^3$ and $X_2X_3X_4$ lie in $\fR$.
Thus $\alpha$  is either $X_2X_3X_5^{i-2}$ or of the form
$X_2X_4^cX_5^{i-c-1}$. Since $p>5$, it follows
from Lemma~\ref{tra}(i) that $X_2X_3X_5^{i-2}$ does not appear in $\tr(X_5^i)$.
Furthermore, since $i\leq p-4$ and $i-c-1\geq 0$, we have $c+3\leq p-2$. Therefore,
by Lemma~\ref{tra}(i), $X_2X_4^cX_5^{i-c-1}$ does not appear in $\tr(X_5^i)$. Thus $\tr(X_5^i)$ does not contribute
to the reduced Gr\"obner basis.

Next assume $i=p-3$. Then the leading monomial of  $\Tr (X_5^i)$ is $X_3^2X_4^{p-5}$. 
Using Lemma~\ref{tra}(i), the only other monomial appearing in $\Tr (X_5^{p-3})$ and not contained in
$\fR$ is $X_2X_4^{p-4}$. The coefficient of $X_3^2X_4^{p-5}$ in $\Tr (X_5^{p-3})$ is
$$\sum_{j\in \bf F_p}{p-3\choose 2}{j\choose 2}^2j^{p-5}={{p-3}\choose 2}\left(\frac{-1}{4}\right)=\frac{-3}{2}$$
and the coefficient of $X_2X_4^{p-4}$ is  
$$\sum_{j\in \bf F_p}(p-3){j\choose 3}j^{p-4}=\frac{-(p-3)}{6}=\frac{1}{2}.$$ 
Therefore 
$$\Tr (X_5^{p-3})\equiv_{\fR} \frac{-3}{2}X_3^2X_4^{p-5}+ \frac{1}{2}X_2X_4^{p-4}.$$ 
Since $X_3^2-X_2(X_3+2X_4)\in \fR$ and $p>5$, it follows that $\fR+(\tr(X_5^{p-3}))=\fR+(X_2X_4^{p-4})$.

For $i=p-2$, the leading monomial  of $\Tr (X_5^i)$ is $X_3X_4^{p-3}$. 
Using  Lemma~\ref{tra}(i), we observe that all monomials less than
$X_3X_4^{p-3}$ which appear in  $\Tr (X_5^{p-2})$
are divisible by at least one of $X_1$, $X_2X_3X_4$, $X_2^2$, $X_3^3$, $X_2X_4^{p-4}$ or $X_3^2X_4^{p-5}$.  
Since all of these monomials are in $\fR+(\tr(X_5^{p-3}))$, it follows that  the contribution to the Hilbert ideal from  
$\Tr (X_5^{p-2})$ is  $X_3X_4^{p-3}$.

For $i=p-1$, the leading monomial of $\Tr (X_5^i)$ is $X_4^{p-1}$. Again 
using  Lemma~\ref{tra}(i),
it is not difficult to see that the smaller monomials appearing in $\tr(X_5^{p-1})$ 
are divisible by at least one of 
$X_1$, $X_2^2$ or $X_3^3$, $X_2X_4^{p-4}$, $X_3^2X_4^{p-5}$ or $X_3X_4^{p-3}$, all of which are in 
$\fR+(\tr(X_5^{p-3}),\tr(X_5^{p-2}))$.
Therefore the contribution to the Hilbert ideal from  $\Tr (X_5^{p-1})$ is $X_4^{p-1}$. 

We define $I:=\fR+(\tr(X_5^i)\mid i =(p-3)/2,\ldots, p-1)$. We have shown that 
$\{X_1,X_2^2, X_3^2-2X_4X_2-X_3X_2, X_4X_3X_2,X_4^{p-4}X_2, X_4^{p-3}X_3, X_4^{p-1}\}$
is a generating set for $I$.
We will show that the remaining families of transfers do not contribute to the Gr\"obner basis and that $N(X_5)$ contributes
$X_5^p$.

\begin{lemma}\label{234lem}
Suppose that $\alpha$ and $\beta$ are monomials with $\alpha<\beta$ and $\deg(\alpha)=\deg(\beta)$. If $X_2X_3X_4$ divides
$\beta$ then $\alpha\in (X_1,X_2^2,X_2X_3^2,X_2X_3X_4)$.
\end{lemma}
\begin{proof}
The lemma follows from the definition of the graded reverse lexicographic order.
\end{proof}

By \cite[3.2 \& 3.6]{shank}, the leading monomial of 
$\tr(X_2X_3X_5^{(p-1)/2})$ is $X_2X_3X_4^{(p-2)/2}$. Therefore, using Lemma~\ref{234lem},
each monomial appearing in $\tr(X_2X_3X_5^{(p-1)/2)})$ lies in $I$. Thus the exceptional transfer
does not contribute to the Gr\"obner basis.

For the fifth family of transfers, using  \cite[3.2 \& 3.6]{shank}, the
leading monomials are $X_2X_3^{p-1-i}X_4^{2i-p-1}$ for $i=(p-1)/2,\ldots,p-2$.
For $i=(p-1)/2$, this gives $X_2X_3^{(p-1)/2}$ which clearly lies in $I$
and, by Lemma~\ref{12}, all of the smaller monomials appearing in $\tr(X_2X_5^{(p-1)/2})$
lie in $I$. For $i>(p-1)/2$, the leading monomial of $\tr(X_2X_5^i)$ is divisible
by $X_2X_3X_4$. Therefore this monomial lies in $I$ and, by Lemma~\ref{234lem},
every monomial appearing in $\tr(X_2X_5^i)$ lies in $I$. Thus the fifth family does
not contribute to the Gr\"obner basis of the Hilbert ideal.
Similarly the invariants of the form $\tr(X_2X_4^2X_5^i)$ appearing in family
three and the invariants of form $\tr(X_2X_4^iX_5^{p-2})$ appearing in family two,
have leading monomials divisible by $X_2X_3X_4$ and therefore do not contribute to the 
Gr\"obner basis.

For the invariants of the form $\Tr (X_4^2X_5^i)$ appearing in family three,
by \cite[3.5]{shank}, the leading monomials are $X_3^{p-1-i}X_4^{2i-p+3}$.
Therefore, as $i$ runs from $(p-1)/2$ to $p-2$, the leading monomials are
$X_3^{(p-1)/2}X_4^2, \ldots, X_3^3X_4^{p-5},X_3^2X_4^{p-3},X_3X_4^{p-1}$.
Clearly these monomials lie in $I$. We will show that the smaller monomials
appearing in these transfers also lie in $I$. Suppose $j>0$ and $\alpha$ is a monomial
with $\alpha<X_3^jX_4^{i+2-j}$, $\deg(\alpha)=i+2$ and $\alpha\not\in I$.
Then one of the following holds:
(i)~$\alpha=X_2X_3X_5^i$,
(ii)~$\alpha=X_2X_4^cX_5^{i+1-c}$ with $c<p-4$,
(iii)~$j=1$, $i=p-2$ and $\alpha=X_3^2X_4^cX_5^{p-2-c}$ with $c<p-3$.
We use Lemma~\ref{tra}(ii). For the first case
$i-d+b+2a=3<p-1$, for the second case $i-d+b+2a=c+1<p-1$
and for the third case $i-d+b+2a=p-2-(p-2-c)+2=c+2<p-1$
Therefore none of these monomials appear in $\tr(X_4^2X_5^i)$.

For the invariants of the form $\tr(X_4^kX_5^{p-2})$ appearing in the
second family, by \cite[3.4]{shank}, the leading monomials are
$X_3X_4^{p+k-3}$ for $k=3,\ldots, p-2$. Clearly these monomials lie in $I$.
We will show that the smaller monomials appearing in these transfers also lie in $I$. 
Suppose $\alpha$ is a monomial
with $\alpha<X_3X_4^{p+k-3}$, $\deg(\alpha)=p+k-2$ and $\alpha\not\in I$.
Then one of the following holds:
(i)~$\alpha=X_2X_3X_5^{p+k-4}$,
(ii)~$\alpha=X_2X_4^cX_5^{p+k-c-3}$ with $c<p-4$,
(iii)~$\alpha=X_3^2X_4^cX_5^{p+k-c-2}$ with $c<p-3$.
Clearly the exponent of $X_5$ must be less than or equal to $p-2$
for any monomial appearing in $\tr(X_4^kX_5^{p-2})$.
For the first case, this exponent on $X_5$ is $p+k-2\geq p+3-2=p-1$.
Thus this monomial does not appear.
Using Lemma~\ref{tra}(ii),
for the second case,
$i-d+b+2a=p-2-(p+k-c-3)+2=c+3-k<p-1$ and for the third case
$i-d+b+2a=p-2-(p+k-c-4)+2=c+4-k<p-1$.
Therefore none of these monomials appear in $\tr(X_4^kX_5^{p-1})$.

For the invariants of the form  $\Tr (X_4^kX_5^{p-1})$ appearing in the first family,
by \cite[3.3]{shank}, the leading monomial is $X_4^{p+k-1}$ for $k=1,\ldots,p-2$.
(The case of $k=0$ appears in the fourth family.)
Clearly these monomials lie in $I$.  As with the previous families,  we will show that all smaller monomials  
appearing in the transfer also lie in $I$. Suppose $\alpha$ is a monomial
with $\alpha<X_4^{p+k-1}$, $\deg(\alpha)=p+k-1$ and $\alpha\not\in I$.
Then one of the following holds:
(i)~$\alpha=X_2X_3X_5^{p+k-2}$,
(ii)~$\alpha=X_2X_4^cX_5^{p+k-c-2}$ with $c<p-4$,
(iii)~$\alpha=X_3^bX_4^cX_5^{p+k-c-1-b}$ with $c<p-3$ and $b=1,2$.
Again we us Lemma~\ref{tra}.
For the first case
$i-d+b+2a=p-1-(p+k-2)+3=4-k<p-1$ and for the second case 
$i-d+b+2a=p-1-(p+k-c-2)+2=c-k+3<p-1$.
Therefore these monomials  do not appear in $\tr(X_4^kX_5^{p-1})$.
For the third case $i-d+b+2a=p-1-(p+k-c-b-1)+b=c-k+2b$. This is less than $p-1$
except for $k=1$, $b=2$, $c=p-4$. However, since $X_3^2-X_2X_3-2X_2X_4\in I$,
we have $X_3^2X_4^{p-4}X_5^2\equiv_I (X_2X_3+2X_2X_4)X_4^{p-4}X_5^2\in I$.
Finally, we consider the invariants of the form  $\Tr (X_2X_4^kX_5^{p-1})$ 
appearing in the first family. Since $\del(X_2)\equiv_{(X_1)} 0$, we have
$\tr(X_2X_4^kX_5^{p-1})\equiv_{(X_1)}X_2 \tr(X_4^kX_5^{p-1})$.
Since we have shown that every monomial appearing in $\tr(X_4^kX_5^{p-1})$
lies in $I$, it follows that every monomial appearing in $\tr(X_2X_4^kX_5^{p-1})$
lies in $I$.

The final element remaining in the generating set is  $N(X_5)$. 
The leading monomial of $N(X_5)$ is clearly $X_5^p$. We will show that the remaining monomials 
appearing in $N(X_5)$ lie in $I$. We choose polynomials $B_0$ and $B_1$ in $\field[X_3,X_4,X_5]$
such that $N(X_5)\equiv_{(X_1,X_2^2)}B_0+X_2B_1$. 
Working modulo $(X_1,X_2)$, the variable $X_5$ generates an $\field\zp$ -- module
isomorphic to $V_3$. Thus we may use the results of Section~\ref{expsec}
to compute $B_0$. By Theorem~\ref{expthm} we have
\begin{eqnarray*}
B_0&\equiv_{(X_3^3)}& X_5^p- X_4^{p-1}X_5+X_3\left(\xi_{11}X_5X_4^{p-2}+\xi_{12}X_5^2X_4^{p-3}\right)\\
&&+X_3^2\left(\xi_{21}X_5X_4^{p-3} +  \xi_{22}X_5^2X_4^{p-4} + \xi_{23}X_5^3X_4^{p-5}\right).
\end{eqnarray*}
Therefore $B_0\equiv_I X_5^p+X_3^2\left( \xi_{22}X_5^2X_4^{p-4} + \xi_{23}X_5^3X_4^{p-5}\right)$.
Since $p>5$, and $X_3^2-X_2X_3-2X_2X_4$ and $X_2X_3X_4$ are both in $I$, we have 
$$X_3^2\left( \xi_{22}X_5^2X_4^{p-4} + \xi_{23}X_5^3X_4^{p-5}\right)\equiv_I
 2\left(\xi_{22}X_5^2X_2X_4^{p-3} + \xi_{23}X_5^3X_2X_4^{p-4}\right).$$
Furthermore, using the fact that $X_2X_4^{p-4}\in I$, gives $B_0\equiv_IX_5^p$.

To complete the proof of Theorem~\ref{v5thm} we need to show that $X_2B_1\in I$.
Note that 
$$N(X_5)\equiv_{(X_1)}\prod_{j\in\fp}\left(X_5+jX_4+{j\choose 2}X_3+{j \choose 3}X_2\right).$$
Therefore
$$X_2B_1=\sum_{j\in\fp}{j \choose 3}X_2\prod_{k\in\fp\setminus\{j\}} \left(X_5+kX_4+{k\choose 2}X_3\right).$$
Since $X_2X_3^2$ and $X_2X_3X_4$ lie in $I$, we have
$$X_2B_1\equiv_I \left(\sum_{j\in\fp}{j \choose 3}X_2\prod_{k\in\fp\setminus\{j\}} \left(X_5+kX_4\right)\right)
+X_2X_3X_5^{p-2}\sum_{j\in\fp}\sum_{k\in\fp\setminus\{j\}}{j\choose 3}{k \choose 2}.$$
Using Lemma~\ref{j4}, we see that $\sum_{k\in\fp\setminus\{j\}}{k \choose 2}=-{j\choose 2}$ giving
$\sum_{j\in\fp}\sum_{k\in\fp\setminus\{j\}}{j\choose 3}{k \choose 2}
=-\sum_{j\in\fp}{j\choose 3}{j \choose 2}=0$ for $p>5$.
Thus
$$X_2B_1\equiv_I \sum_{j\in\fp}{j \choose 3}X_2\prod_{k\in\fp\setminus\{j\}} \left(X_5+kX_4\right).$$
However
\begin{eqnarray*}
 \prod_{k\in\fp\setminus\{j\}} \left(X_5+kX_4\right)&=&\frac{\prod_{k\in\fp} \left(X_5+kX_4\right)}{X_5+jX_4}\\
&=&\frac{X_5^p-X_5X_4^{p-1}}{X_5+jX_4} = X_5^{p-1}\left(\frac{1-(X_4/X_5)^{p-1}}{1+(jX_4/X_5)}\right).
\end{eqnarray*}
For the purposes of computing $X_2B_1$ modulo $I$, we may assume $j\not=0$. This means that $(j)^{p-1}=1$. Thus,
using $(1-a^n)/(1-a)=1+a+ \dots + a^{n-1}$ with $a=-jX_4/X_5$,
we see that  
\begin{eqnarray*}
 \prod_{k\in\fp\setminus\{j\}} \left(X_5+kX_4\right)
&=&X_5^{p-1}\left(1+\left(-jX_4/X_5\right)+\cdots +\left(-jX_4/X_5\right)^{p-2}\right)\\
&=&X_5^{p-1}-jX_4X_5^{p-2}+j^2X_4^2X_5^{p-3}+\cdots +(-j)^{p-2}X_4^{p-2}X_5.
\end{eqnarray*}
Therefore
$$X_2B_1\equiv_I \sum_{j\in\fp}
{j \choose 3}X_2\left(X_5^{p-1}-jX_4X_5^{p-2}+\cdots +(-j)^{p-2}X_4^{p-2}X_5\right).$$
Since ${j \choose 3}$ is a polynomial of degree $3$ in $j$, using Lemma~\ref{j4} gives
$$X_2B_1\equiv_I \sum_{j\in\fp}
{j \choose 3}X_2\left((-j)^{p-4}X_4^{p-4}X_5^3+(-j)^{p-3}X_4^{p-3}X_5^2+(-j)^{p-2}X_4^{p-2}X_5\right).$$
Therefore, since $X_2X_4^{p-4}\in I$, we have $X_2B_1\in I$.

We have shown that $N(X_5)\equiv_I X_5^p$. Therefore the Hilbert ideal is generated by
$$\{X_1,X_2^2, X_3^2-2X_4X_2-X_3X_2, X_4X_3X_2,X_4^{p-4}X_2, X_4^{p-3}X_3, X_4^{p-1}, X_5^p\}.$$ 
It is clear that this set is a reduced Gr\"obner basis.
The corresponding monomial basis consists of all monomials not divisible by any of the
generators and the description of the Hilbert series comes from the monomial basis.

\begin{remark}
We observe that the top degree of $\field[V_5]_{\zp}$ is $2p-3$.
It is clear that $2p-3$ is an upper bound for the Noether number of
$V_4$. It follows from Remark~\ref{v4rem} and \cite[4.2]{sw:noeth},
that the Noether number of $V_5$ is $2p-3$.
\end{remark}

\section{The module structure  for the coinvariants of $V_4$ and $V_5$}

In this section we use the bases  constructed in Sections~\ref{v4coin} and \ref{v5coin}
to determine the $\field\zp$ -- module structure of the coinvariants of $V_4$ and $V_5$.
Note that, since the Hilbert ideal is homogeneous, the coinvariants are a graded ring. Furthermore,
the group action preserves degrees. Thus the homogeneous components are $\field\zp$ -- module summands.
We will refine this decomposition by describing each homogeneous component as a direct sum of
indecomposable modules. Recall that the socle of a module is the sum of its irreducible submodules.
For an $\field\zp$ -- module, this is the span of the fixed points. A non-zero cyclic $\field\zp$ -- module
has a one dimensional socle and,  since all indecomposable $\field\zp$ -- modules are cyclic,
the dimension of the socle is the number of summands. For a non-zero cyclic module with socle
${\rm Span}(v)$, we will say that $v$ {\it determines} the socle.

\begin{lemma}\label{soclelem}
Suppose that $W_1, W_2,\ldots W_m$ are cyclic submodules of $W$ and that
$\omega_i$ determines the socle of $W_i$.
If $\{\omega_1,\omega_2,\ldots,\omega_m\}$ is linearly independent
and $\dim(W)=\dim(W_1)+\dim(W_2)+\cdots +\dim(W_m)$, then
$W=W_1\oplus W_2\oplus \cdots \oplus W_m$.
\end{lemma}
\begin{proof}
For a homomorphism of modules,
the socle of the kernel is the kernel of the restriction of the homomorphism to the socle.
Thus a homomorphism which is injective on its socle is injective.
Apply this to the homomorphism from the external direct sum of the $W_i$ to their internal sum.
Since $\{\omega_1,\omega_2,\ldots,\omega_m\}$ is linearly independent, this map is
injective on its socle and hence injective. Therefore the internal sum of the $W_i$ is direct and
$W_1\oplus W_2\oplus \cdots \oplus W_m$ is a subspace of $W$. However,
since $\dim(W)=\dim(W_1)+\dim(W_2)+\cdots \dim(W_m)$, the subspace coincides with $W$.
\end{proof}
  
We define the {\it weight} of a monomial in $\field[V_n]$ by
$\wt(X_1^{e_1}\cdots X_n^{e_n})=e_1+2e_2+ \dots +ne_n.$ 
If $f$ is a linear combination of monomials of the same weight,
we will refer to $f$ as {\it isobaric} and we will take the weight of $f$ to be the
common weight of the monomials appearing in $f$.
Note that
if $\beta$ is a monomial appearing
in $\Delta (f)$ with $f$ isobaric, then
$\wt(\beta)<\wt(f)$. Thus, for a fixed positive integer $m$, the span of the monomials of weight
less than $m$ forms an $\field\zp$ -- submodule. Allowing $m$ to vary over the positive integers gives
a weight filtration of the polynomial ring. 
For $V_4$ and $V_5$ we fix a basis for the coinvariants given by images of
monomials. For $V_4$, the basis is given in Theorem~\ref{h4} and for $V_5$ the basis is given by
Theorem~\ref{v5thm}. We define the weight of the basis elements to be the weight of the corresponding monomial
and, as in the polynomial ring, a linear combination of basis elements of a common weight is isobaric
with a well defined weight.

\begin{lemma}\label{coinwtlem} If $f$ is an isobaric coinvariant of weight $m$, then $\del(f)$ is in the span of the
basis elements of weight less than $m$.
\end{lemma}
\begin{proof}
Since $\del$ is linear it is sufficient to consider $\del(\beta)$ for a basis element $\beta$ of weight $m$.
To compute $\del(\beta)$, we lift to the corresponding monomial in the polynomial ring, say $\ov{\beta}$, compute $\del(\ov{\beta})$,
and then project back to coinvariants. The terms appearing in $\del(\ov{\beta})$ all have weight less than $m$.
For $V_4$, the reduced Gr\"obner basis is a set of monomials.
Thus each term appearing in $\del(\ov{\beta})$ either projects to zero or
projects to a term of weight less than $m$. For $V_5$, there are seven monomial relations and
one non-isobaric relation given by
$X_3^2-2X_2X_4-X_2X_3$. This last relation is used to give a rewriting rule which replaces the product $x_3\cdot x_3$
with $2x_2x_4+x_2x_3$. Thus an element of weight $6$ in the polynomial ring is identified with a sum of
two terms, one of weight $6$ and one of weight $5$, in the coinvariants. Thus each term appearing in $\del(\ov{\beta})$
either projects to zero or projects to a linear combination of terms with weight less than $m$. 
\end{proof}

As a consequence of Lemma~\ref{coinwtlem}, for each positive integer $m$, the span of the basis elements of weight less
than $m$ form an $\field\zp$ -- submodule. Collectively these submodules give a weight filtration of the coinvariants.
Suppose $\beta$  is a basis element of weight $m$. Define $\delta(\beta)$ to be the sum of terms of weight $m-1$ appearing in $\del(\beta)$
and extend $\delta$ to linear map on the coinvariants. We can think of $\delta$ as the linear map induced by $\del$ on 
the associated graded module of the weight filtration. In the following we use $\field[V]_{\zp}^d$ to denote the
homogeneous component of degree~$d$.

\begin{lemma}\label{wtlem} Suppose $n$ is $4$ or $5$, and
$m$ is the minimum weight occurring in $\field[V_n]_{\zp}^d$.
For an isobaric coinvariant $f$ of weight $\ell$ and a positive integer $k$,
 any term appearing in $\delta^k(f)-\del^k(f)$ has weight less than $\ell-k$. In particular,
if $\ell=m+k$, then $\delta^k(f)=\del^k(f)$. Furthermore, if $\ell=m$, then $f$ is invariant.
\end{lemma}
\begin{proof}
The proof is by induction on $k$. For $k=1$, the result is essentially the definition of $\delta$.
Suppose the result is true for $k>1$. Then $\delta^k(f)=\del^k(f)+h$ where $h$ is  a sum of terms of
weight less than $\ell-k$. Thus $\delta(\delta^k(f))$ consists of the sum of the terms of weight
$\ell-k-1$ in $\del(\del^k(f))+\del(h)$. However, from Lemma~\ref{coinwtlem}, all of terms appearing
in $\del(h)$ have weight less than $\ell-k-1$. Therefore $\delta^{k+1}(f)$ consists of the sum of the
terms of weight $\ell-(k+1)$ appearing in $\del^{k+1}(f)$, as required.
If $\ell-k=m$, there are no terms of weight less than $\ell-k$ so  $\delta^k(f)=\del^k(f)$.
If $\ell=m$, the fact that $f$ is invariant follows from Lemma~\ref{coinwtlem}.
\end{proof}

The following lemma will play an important role in determining the $\field\zp$ -- module structure of
$\ivcoin$.

\begin{lemma}\label{4vdelta} In $\ivcoin$, for $j\geq k$,
$$\delta^k(x_3^ix_4^j)=\frac{j!}{(j-k)!}x_3^{i+k}x_4^{j-k}+\frac{j!}{(j-k+1)!}\left(i k +{k \choose 2}\right)x_2 x_3^{i+k-2} x_4^{j-k+1}.$$
\end{lemma}
\begin{proof} The proof is by induction on $k$. For $k=1$, a straight forward calculation gives
$\delta(x_3^ix_4^j)=j x_3^{i+1} x_4^{j-1} + i x_2 x_3^{i-1} x_4^j$.
For $k\geq 1$ we have
\begin{eqnarray*}
 \delta^{k+1}(x_3^ix_4^j) &=&\delta(\delta^k(x_3^ix_4^j))\\
 &=& \delta\left(\frac{j!}{(j-k)!}x_3^{i+k}x_4^{j-k}+\frac{j!}{(j-k+1)!}\left(i k +{k \choose 2}\right)x_2 x_3^{i+k-2} x_4^{j-k+1}\right)\\
 &=&\frac{j!}{(j-k)!}(j-k)x_3^{i+k+1}x_4^{j-k-1} +\frac{j!}{(j-k)!}(i+k)x_2 x_3^{i+k-1} x_4^{j-k}\\
 & &+\frac{j!}{(j-k+1)!}\left(i k +{k \choose 2}\right)(j-k+1)x_2 x_3^{i+k-1} x_4^{j-k}\\
 &=& \frac{(j!) x_3^{i+(k+1)}x_4^{j-(k+1)}}{(j-(k+1))!}+\frac{(j!) x_2 x_3^{i+k-1} x_4^{j-k}}{(j-k)!}\left((i+k)+ik +{k \choose 2}\right)\\
 &=&\frac{(j!) x_3^{i+(k+1)}x_4^{j-(k+1)}}{(j-(k+1))!} 
+\frac{(j!) x_2 x_3^{i+(k+1)-2} x_4^{j-(k+1)+1}}{(j-(k+1)+1)!}\left(i(k+1)+{k+1 \choose 2}\right),
\end{eqnarray*}

\noindent as required.
\end{proof}

\begin{theorem} (i) $\ivcoin^0\cong\ivcoin^{2p-3}\cong V_1$, $\ivcoin^1\cong V_3$.

\noindent (ii) For $d=p,\ldots, 2p-4$,
$$\ivcoin^d=x_3^{d-(p-1)}x_4^{p-1}\field\zp\oplus\ x_2x_3^{d-p}x_4^{p-1}\field\zp\cong V_{2p-2-d}\oplus V_{2p-3-d}$$
with $\left(\ivcoin^d\right)^{\zp}={\rm Span}\{x_3^{p-2}x_4^{d-(p-2)}, x_2x_3^{p-4}x_4^{d-(p-3)}\}$.

\noindent (iii) For $d=p-1, p-2$,
$$\ivcoin^d=x_4^d\field\zp\oplus\ x_2x_4^{d-1}\field\zp\cong V_{p-1}\oplus V_{p-3}$$
with $\left(\ivcoin^d\right)^{\zp}={\rm Span}\{x_3^{p-2}x_4^{d-(p-2)}, x_2x_3^{p-4}x_4^{d-(p-3)}\}$.

\noindent (iv) For $d=2,\ldots,p-3$,
$$\ivcoin^d=x_4^d\field\zp\oplus\ \left(x_3^2x_4^{d-2}-\frac{d+2}{2} x_2x_4^{d-1}\right)\field\zp\cong V_{d+2}\oplus V_{d-1}$$
with $\left(\ivcoin^d\right)^{\zp}={\rm Span}\{x_3^d-d x_2x_3^{d-2}x_4,x_2 x_3^{d-1}\}$.
\end{theorem}

\begin{proof} Part (i) is clear.

(ii) For $p\leq d\leq 2p-4$, from Theorem~\ref{h4}, a basis for $\ivcoin^d$ is given by
$$x_3^{p-2}x_4^{d-(p-2)},\ x_3^{p-3}x_4^{d-(p-3)},\ldots,x_3^{d-(p-1)}x_4^{p-1}$$ and
$$x_2x_3^{p-4}x_4^{d-(p-3)},\ x_2x_3^{p-5}x_4^{d-(p-4)}, \ldots,x_2x_3^{d-p}x_4^{p-1}.$$ Therefore
the dimension of $\ivcoin^d$ is $(2p-2-d)+(2p-3-d)$. The elements $x_3^{p-2}x_4^{d-(p-2)}$ and $x_2x_3^{p-4}x_4^{d-(p-3)}$
are invariant and have minimum weight. From Lemma~\ref{4vdelta}, 
$\delta^{2p-3-d} (x_3^{d-(p-1)}x_4^{p-1})$ is a linear combination of
$x_3^{p-2}x_4^{d-(p-2)}$ and $x_2x_3^{p-4}x_4^{d-(p-3)}$ with the
coefficient of $x_3^{p-2}x_4^{d-(p-2)}$ non-zero. Applying Lemma~\ref{wtlem} gives
$\Delta^{2p-3-d} (x_3^{d-(p-1)}x_4^{p-1})=\delta^{2p-3-d} (x_3^{d-(p-1)}x_4^{p-1})$.
Thus $x_3^{p-2}x_4^{d-(p-2)}$
generates a module of dimension $2p-2-d$. Again using Lemma~\ref{4vdelta},
$$\delta^{2p-4-d}(x_2x_3^{d-p}x_4^{p-1})=x_2 \delta^{2p-4-d}(x_3^{d-p}x_4^{p-1})=cx_2x_3^{p-4}x_4^{d-(p-3)}$$
with $c=(p-1)!/(d-p+3)!\not=0$. Therefore $x_2x_3^{d-p}x_4^{p-1}$ generates a module of dimension
$2p-3-d$. Since the fixed points are linearly independent, using Lemma~\ref{soclelem} shows that
the sum of $x_3^{d-(p-1)}x_4^{p-1}\field\zp$ and $\ x_2x_3^{d-p}x_4^{p-1}\field\zp$ is direct.
Therefore $x_3^{d-(p-1)}x_4^{p-1}\field\zp+\ x_2x_3^{d-p}x_4^{p-1}\field\zp$ is a submodule 
isomorphic to $V_{2p-2-d}\oplus V_{2p-3-d}$. Since the dimensions match, this submodule is 
all of $\ivcoin^d$.

(iii) The proof for case (iii) is similar to case (ii). It follows
from Lemma~\ref{4vdelta} that $x_4^d$ and $x_2x_4^{d-1}$ generate modules of the
dimensions $p-1$ and $p-3$, respectively and that $\delta^{p-2}(x_4^d)$ and $\delta^{p-4}(x_2x_4^{d-1})$
are linearly independent invariants of minimum weight. Therefore, using Lemma~\ref{soclelem}, we have identified a submodule
isomorphic to $V_{p-1}\oplus V_{p-3}$. The result follows from the observation that both
$\ivcoin^{p-1}$ and $\ivcoin^{p-2}$ have dimension $2p-4$.

(iv) For $2\leq d\leq p-3$, a basis for $\ivcoin^d$ is given by
 $x_3^d,x_3^{d-1}x_4,\ldots,x_4^d$ and
$x_2x_3^{d-1}$, $x_2x_3^{d-2}x_4,\ldots,x_2x_4^{d-1}$. Therefore
the dimension of $\ivcoin^d$ is $2d+1$. The minimum weight subspace is given by
${\rm Span}(x_2x_3^{d-1})$. A second invariant, isobaric but with  non-minimum weight, is given by 
$x_3^d-dx_2x_3^{d-2}x_4$.
Using Lemma~\ref{4vdelta}, $\delta^d(x_4^d)=d!\left(x_3^d+d(d-1)x_2x_3^{d-2}x_4/2\right)$.
Direct calculation gives 
$\delta^{d+1}(x_4^d)=\delta\left(\delta^d(x_4^d)\right)=d!d(d-1)x_2x_3^{d-1}/2$.
Thus $x_4^d$ generates a module of dimension $d+2$. 
Again using Lemma~\ref{4vdelta} gives
$\delta^{d-2}(x_4^d-(d+2) x_2x_4^{d-1}/2)=(d-1)!(x_3^d-d x_2x_3^{d-2}x_4)$.
The only basis element with weight less than $3d$ is the invariant $x_2x_3^{d-1}$. Therefore, using Lemma~\ref{wtlem},
$\Delta^{d-2}(x_4^d-(d+2) x_2x_4^{d-1}/2)=(d-1)!(x_3^d-d x_2x_3^{d-2}x_4)+cx_2x_3^{d-1}$
for some constant $c$.
Thus $x_4^d-(d+2) x_2x_4^{d-1}/2$
generates a module of dimension $d-1$. The intersection of the socles of
the two given submodules is trivial. Therefore, using Lemma~\ref{soclelem}, the sum of the modules is direct.
Thus $\ivcoin^d$ has a submodule isomorphic to $V_{d+2}\oplus V_{d-1}$ and,
since the dimension of $\ivcoin^d$ is $2d+1$, this submodule is all of
$\ivcoin^d$.
\end{proof}

We will require a number of technical lemmas to determine the $\field\zp$ -- module structure of
$\vcoin$.


\begin{lemma}\label{345dellem} In $\vcoin$, for $j\geq k$, 
$$\delta^k(x_3x_4^ix_5^j)=\frac{j!}{(j-k)!}x_3 x_4^{k+i} x_5^{j-k} + \frac{j!(2ik+k^2)}{(j-k+1)!}x_2x_4^{k+i-1}x_5^{j-k+1} +c_k x_2x_3x_5^{j-k+2}$$
where $c_k=0$ unless $k+i=3$ in which case $c_k$ equals the coefficient of $x_2x_4^{i+k-2}x_5^{j-k+2}$ in $\delta^{k-1}(x_3x_4^ix_5^j)$.
\end{lemma}
\begin{proof}
The proof is by induction on $k$. First consider $k=1$. A direct calculation gives
$\delta(x_3 x_4^i x_5^j)=x_2 x_4^i x_5^j+i x_3^2 x_4^{i-1}+ j x_3 x_4^{i+1}  x_5^{j-1}$.
In $\vcoin$, we have the  relation $x_3^2=x_2(2 x_4+ x_3)$. Since $\delta$ picks out the highest weight terms of $\del$
we may substitute $2x_2 x_4$ for $x_3^2$ giving
$\delta(x_3 x_4^i x_5^j)=(2i+1)x_2 x_4^i x_5^j+ j x_3 x_4^{i+1}  x_5^{j-1}$.
For $k\geq 1$, we have 
\begin{eqnarray*}
\delta^{k+1}(x_3x_4^ix_5^j)&=&\delta\left(\delta^k(x_3x_4^ix_5^j)\right)\\
 &=&\delta\left(\frac{j!}{(j-k)!}x_3 x_4^{k+i} x_5^{j-k} + \frac{j!(2ik+k^2)}{(j-k+1)!}x_2x_4^{k+i-1}x_5^{j-k+1} +c_k x_2x_3x_5^{j-k+2}\right)\\
 &=& \frac{j!}{(j-k-1)!}x_3 x_4^{k+i+1} x_5^{j-k-1}+\frac{j!(k+i)}{(j-k)!}x_3^2 x_4^{k+i-1} x_5^{j-k}+ \frac{(j!)x_2 x_4^{k+i} x_5^{j-k}}{(j-k)!}\\
 && + \frac{j!(2ik+k^2)}{(j-k)!}x_2x_4^{k+i}x_5^{j-k} + c_{k+1} x_2x_3x_5^{j-k+2}.
\end{eqnarray*}
Substituting $2x_2 x_4$ for $x_3^2$ gives
\begin{eqnarray*}
\delta^{k+1}(x_3x_5^j)&\equiv_{(x_2x_3)}&\frac{(j!)x_3 x_4^{k+i+1} x_5^{j-k-1}}{(j-k-1)!}+\frac{j!(2(k+i)+1+2ik+k^2)}{(j-k)!}x_2 x_4^{k+i} x_5^{j-k}\\
&\equiv_{(x_2x_3)}&\frac{(j!)x_3 x_4^{i+k+1} x_5^{j-(k+1)}}{(j-k-1)!}+\frac{j!(2i(k+1)+(k+1)^2)}{(j-k)!}x_2 x_4^{k+i} x_5^{j-k}
\end{eqnarray*}
as required.
\end{proof}

\begin{lemma}\label{35dellem}
For $p-4\geq d>3$, $\delta^d(x_3x_5^{d-1})=d(d!)x_2 x_4^{d-1}$.
\end{lemma}
\begin{proof} From Lemma~\ref{345dellem}, 
$$\delta^{d-1}(x_3 x_5^{d-1})\equiv_{(x_2 x_3)}(d-1)!x_3 x_4^{d-1}+(d-1)!(d-1)^2 x_2 x_4^{d-2}x_5.$$
Applying $\delta$ and using Lemma~\ref{345dellem} gives
$$\delta^d(x_3 x_5^{d-1})=(d-1)!\left((2d-1)+(d-1)^2\right)x_2x_4^{d-1}=d^2(d-1)!x_2x_4^{d-1}$$
as required.
\end{proof}

\begin{lemma}\label{45dellem} In $\vcoin$, for $j\geq k$, 
$$\delta^k(x_4^i x_5^j) = a_k x_4^{i+k} x_5^{j-k} + b_k x_3 x_4^{i+k-2} x_5^{j-k+1}
   + c_k x_2 x_4^{i+k-3}x_5^{j-k+2} +d_k x_2 x_3 x_5^{i+j-2}$$
where
$$a_k=\frac{j!}{(j-k)!},\hspace{1cm}  b_k=\frac{j!}{(j-k+1)!}\left(i k +{k \choose 2}\right),$$

$$c_k=\frac{j!}{(j-k+2)!}{k \choose 2}\left(2i^2+ (2k-5)i + \frac{(k-2)(3k-7)}{6}\right)$$
and $d_k=0$ unless $i+k=5$ in which case $d_k=c_{k-1}$.
\end{lemma}
\begin{proof}
 The proof is by induction on $k$. For $k=1$, a straight forward calculation gives
$\delta(x_4^ix_5^j)=j x_4^{i+1} x_5^{j-1} + i x_3 x_4^{i-1} x_5^j$.
For $k\geq 1$ we have 
$$\delta^{k+1}(x_4^i x_5^j) = \delta\left(a_k x_4^{i+k} x_5^{j-k} + b_k x_3 x_4^{i+k-2} x_5^{j-k+1}
   + c_k x_2 x_4^{i+k-3}x_5^{j-k+2} +d_k x_2 x_3 x_5^{i+j-2}\right).$$
Using the definition of $\delta$ and Lemma~\ref{345dellem} gives
\begin{eqnarray*}
\delta^{k+1}(x_4^i x_5^j) 
 &=& a_k(j-k) x_4^{i+k+1} x_5^{j-k-1} + a_k (i+k) x_3 x_4^{i+k-1} x_5^{j-k}\\
 && + b_k(j-k+1)x_3 x_4^{i+k-1} x_5^{j-k} + b_k (2(i+k-2)+1) x_2 x_4^{i+k-2} x_5^{j-k+1}\\
 && + c_k (i+k-3) x_2 x_3 x_4^{i+k-4} x_5^{j-k+2} + c_k (j-k+2) x_2 x_4^{i+k-2} x_5^{j-k+1}\\
 &=&  a_k(j-k) x_4^{i+k+1} x_5^{j-k-1} + \left(a_k(i+k)+b_k(j-k+1)\right)x_3 x_4^{i+k-1} x_5^{j-k}\\
 && + \left(b_k(2i+2k-3)+c_k(j-k+2)\right) x_2 x_4^{i+k-2} x_5^{j-k+1}\\
 && + c_k (i+k-3) x_2 x_3 x_4^{i+k-4} x_5^{j-k+2}.
\end{eqnarray*}
Since $x_2x_3x_4=0$, it is clear that $d_{k+1}=0$ unless $i+(k+1)-5=0$ in which case $d_{k+1}=c_k$.
The fact that $a_{k+1}= a_k(j-k)$ and $b_{k+1}=a_k(i+k)+b_k(j-k+1)$ follows from the proof of
Lemma~\ref{4vdelta}. Thus we need only verify $c_{k+1}$. The coefficient of 
$x_2 x_4^{i+(k+1)-3} x_5^{j-(k+1)+2}$ in the preceding expression is
$b_k(2i+2k-3)+c_k(j-k+2)$. Substituting the expressions
for $b_k$ and $c_k$ gives
{\small
$$
 \frac{(j!)(2i+2k-3)}{(j-k+1)!}\left(i k +{k \choose 2}\right)
+\frac{(j!)(j-k+2)}{(j-k+2)!}{k \choose 2}\left(2i^2+ (2k-5)i + \frac{(k-2)(3k-7)}{6}\right).
$$}
Factoring gives
{\small
$$
 \frac{(j!)}{(j-k+1)!}\left( 
(2i+2k-3)\left(i k +{k \choose 2}\right)
+{k \choose 2}\left(2i^2+ (2k-5)i + \frac{(k-2)(3k-7)}{6}\right)\right).
$$}
A MAGMA\cite{magma} calculation can be used to verify that, as polynomials in $i$ and $k$,
$$
{k+1 \choose 2}\left(2i^2+(2(k+1)-5)i+\frac{\left((k+1)-2\right)\left(3(k+1)-7\right)}{6}\right)$$
equals
$$
(2i+2k-3)\left(i k +{k \choose 2}\right)
+{k \choose 2}\left(2i^2+ (2k-5)i + \frac{(k-2)(3k-7)}{6}\right).
$$
This completes the induction step.
\end{proof}

\begin{lemma}\label{5dellem}
(i) For $p-3\geq d>4$, 
$$\delta^{d+1}(x_5^d)=\frac{d(d+1)!}{12}\left(6x_3x_4^{d-1}+(d-1)(3d-4)x_2 x_4^{d-2}x_5\right).$$

\noindent (ii) For $p-4\geq d>3$, $\delta^{d+2}(x_5^d)=\frac{d(3d-1)(d+2)!}{12}x_2 x_4^{d-1}$.
\end{lemma}
\begin{proof} 
Using Lemma~\ref{45dellem}, 
$$ \delta^d(x_5^d)\equiv_{(x_2 x_3)} d!x_4^d+d!{d \choose 2}\left( x_3 x_4^{d-2} x_5+\frac{(d-2)(3d-7)}{12} x_2x_4^{d-3}x_5^2\right).$$
Applying $\delta$ and using Lemma~\ref{345dellem} gives
\begin{eqnarray*}
\frac{\delta^{d+1}(x_5^d)}{d!}&\equiv_{(x_2 x_3)}&\left(d+{d \choose 2}\right)x_3x_4^{d-1}
+{d \choose 2}\left(2d-3 + \frac{3d^2 -13 d + 14}{6}\right)x_2 x_4^{d-2}x_5\\
&\equiv_{(x_2 x_3)}&\frac{d(d+1)}{2}x_3x_4^{d-1}
+{d \choose 2}\frac{3d^2 - d - 4}{6}x_2 x_4^{d-2}x_5\\
&\equiv_{(x_2 x_3)}&\frac{d(d+1)}{12}\left(6x_3x_4^{d-1}+(d-1)(3d-4)x_2 x_4^{d-2}x_5\right)
\end{eqnarray*}
Therefore,
$$\delta^{d+1}(x_5^d)=\frac{d(d+1)(d!)}{12}\left(6x_3x_4^{d-1}+(d-1)(3d-4)x_2 x_4^{d-2}x_5\right)+c x_2x_3x_5^{d-2}$$
where $c=0$ unless $d=4$ in which case $c=(d!)d(d-1)(d-2)(3d-7)/24=120$. Again applying $\delta$
and using Lemma~\ref{345dellem} gives
\begin{eqnarray*}
\delta^{d+2}(x_5^d)&\equiv_{(x_2 x_3)}&\frac{d(d+1)(d!)}{12}\left(6(2(d-1)+1)+(d-1)(3d-4)\right)x_2 x_4^{d-1}\\
 &\equiv_{(x_2 x_3)}&\frac{d(d+1)(d!)(3d^2+5d-2)}{12}x_2 x_4^{d-1}.
\end{eqnarray*}
Therefore,
$$\delta^{d+2}(x_5^d)=\frac{d(d+1)(d!)(3d-1)(d+2))}{12}x_2 x_4^{d-1}+c'x_2x_3x_4^{d-2}$$
where $c'=0$ unless $d=3$ in which case $c'=(d!)d(d-1)(d+1)(3d-4)/12=60$.
\end{proof}

\begin{lemma}\label{numlem}
For $1<d<p$, $\vcoin^d$ is generated as an $\field\zp$ -- module by
$$\{x_5^d, x_3x_5^{d-1}, x_2x_5^{d-1}, x_2x_3x_5^{d-2}\}.$$
Hence $\vcoin^d$ decomposes into a sum of at most four indecomposable summands.
\end{lemma}
\begin{proof} Having fixed a basis for $\vcoin$ consisting of the images of monomials
we can use the order on $\field[V_5]$ to give a total order on the basis and a partial order
on the coinvariants. Thus it is possible to determine the leading term of a coinvariant.
Note, however, that the order is not multiplicative. We will denote the leading term of a coinvariant $f$
by $\lt(f)$.
To show that  $\vcoin^d$ is generated by
$\Gamma:=\{x_5^d, x_3x_5^{d-1}, x_2x_5^{d-1}, x_2x_3x_5^{d-2}\}$
it is sufficient to show that 
$\{\lt(\Delta^k(\beta))\mid \beta \in \Gamma \}$ spans $\vcoin^d$.
Furthermore, for every $\beta\in\Gamma$, $\lt(\delta^k(\beta))=\lt(\Delta^k(\beta))$.
Observe that $\delta^k(x_2x_5^{d-1})=x_2\delta^k(x_5^{d-1})$. Thus 
$\delta^k(x_2x_5^{d-1})$ can be computed for $k\leq d-1$  using Lemma~\ref{45dellem}.
Therefore, using  Lemma~\ref{45dellem} and Lemma~\ref{345dellem}, we see that
$\{ x_2x_3x_5^{d-2}\} \cup \{\lt(\delta^k(x_5^d))\mid 0\leq k\leq d\}
\cup \{\lt(\delta^k( x_3x_5^{d-1})),\lt(\delta^k(x_2x_5^{d-1}))\mid 0\leq k\leq d-1\}$
is a basis for $\vcoin^d$. Hence $\Gamma$ is a generating set.

To see that the number of generators is an upper bound on the number of indecomposable summands,
work inductively. Certainly a module with one generator is indecomposable. 
Suppose a module has more than one generator. It is convenient to define the length of a generator
to be the dimension of the submodule it generates.
By looking at the decomposition of the module, it is not hard to see that a generator of maximum length
generates a summand.
\end{proof}

\begin{theorem} Suppose $p>5$.

\noindent (i) $\vcoin^0\cong\vcoin^{2p-3}\cong V_1$ and $\vcoin^1\cong V_4$.

\noindent (ii) $\vcoin^{2p-4}\cong V_2\oplus V_1$ and $\vcoin^2\cong V_6\oplus V_2$.

\noindent (iii) For $d=p+2,\ldots, 2p-5$:  $\vcoin^d\cong V_{2p-d-2} \oplus V_{2p-d-3}\oplus V_{2p-d-4}$.

\noindent (iv) For $d=p,\,p+1$: $\vcoin^d\cong V_{2p-d-2} \oplus V_{2p-d-3}\oplus V_{2p-d-4}\oplus V_1$.

\noindent (v) For $d=p-1,\,p-2$: $\vcoin^d\cong V_{p-1} \oplus V_{p-3}\oplus V_{p-4}\oplus V_1$.

\noindent (vi) For $d=p-3$ and $p>11$: $\vcoin^d\cong V_{p-1} \oplus V_{p-3}\oplus V_{p-5}\oplus V_1$.

\noindent (vii) For $d=5,\ldots,p-4$: if $3d-1\not\equiv_{(p)} 0$ and $3d-2\not\equiv_{(p)}0$ then
$$\vcoin^d\cong V_{d+3} \oplus V_d\oplus V_{d-2}\oplus V_1;$$
if $3d-1\not\equiv_{(p)} 0$ and $3d-2\equiv_{(p)}0$ then
$$\vcoin^d\cong V_{d+3} \oplus 2 V_{d-1} \oplus V_1;$$
if $3d-1\equiv_{(p)} 0$ then
$\vcoin^d\cong V_{d+2} \oplus V_{d+1}\oplus V_{d-2}\oplus V_1$.

\noindent (viii) For $p>11$:  $\vcoin^3 \cong V_6\oplus V_4 \oplus V_1$ and $\vcoin^4\cong V_7\oplus V_4 \oplus V_3$.
\end{theorem}

\begin{remark} MAGMA \cite{magma} calculations give the following.

\noindent (i) For $p=5$, the homogeneous component of $\field[V_5]_{{\bf Z}/5}$ in increasing degree are
isomorphic to
 $V_1$, $V_4$, $2V_4$, $2V_4\oplus 2V_1$, $2V_4\oplus 2V_1$, $V_3\oplus V_4\oplus2V_1$, $V_4\oplus 2V_1$, $2V_1$.

\noindent (ii) For $p=11$: 
$\field[V_5]_{{\bf Z}/11}^3 \cong V_6\oplus V_4 \oplus V_1$, $\field[V_5]_{{\bf Z}/11}^4\cong V_6\oplus V_5 \oplus V_3$ and
 $\field[V_5]_{{\bf Z}/11}^8\cong V_{10}\oplus 2V_7 \oplus V_1$

\noindent (iii) For $p=7$: $\field[V_5]_{{\bf Z}/7}^3 \cong V_6\oplus V_3 \oplus V_2$ and 
$\field[V_5]_{{\bf Z}/7}^4\cong V_6\oplus V_4 \oplus V_3$.
\end{remark}

\begin{proof} Part (i) is clear.

(ii) For $d=2p-4$: $x_4^{p-3} x_5^{p-1}$ generates a submodule of dimension $2$ and both $x_4^{p-2} x_5^{p-2}$ and 
$x_3x_4^{p-4}x_5^{p-1}$ are invariant. For $d=2$: A straight forward calculation shows that $x_5^2$ generates a submodule
of dimension $6$ with socle ${\rm Span}(x_2x_3)$. A second calculation shows that $2x_4^2-3x_3x_5-3x_2x_5-2x_2x_4$ generates a submodule
of dimension $2$ with socle ${\rm Span}(x_3x_4-3x_2 x_5-2x_2x_4)$. Since the dimension of the degree $2$ homogeneous component
is $8$, Lemma~\ref{soclelem} applies to give the stated decomposition.

(iii) For $p+2\leq d \leq 2p-5$, a basis for $\vcoin^d$ is given by 
\begin{eqnarray*}
&&x_4^{p-2}x_5^{d-p+2}, x_4^{p-3}x_5^{d-p+3}, \dots, x_4^{d-p+1}x_5^{p-1},\\
&&x_3x_4^{p-4}x_5^{d-p+3}, x_3x_4^{p-5}x_5^{d-p+4}, \dots ,x_3x_4^{d-p}x_5^{p-1},\\
&&x_2x_4^{p-5}x_5^{d-p+4}, x_2x_4^{p-6}x_5^{d-p+5}, \dots ,x_2x_4^{d-p}x_5^{p-1}.
\end{eqnarray*}
Therefore the dimension of $\vcoin^d$ is $(2p-d-2)+(2p-d-3)+(2p-d-4)$. The elements  $x_4^{p-2}x_5^{d-p+2}$,
$x_3x_4^{p-4}x_5^{d-p+3}$ and $x_2x_4^{p-5}x_5^{d-p+4}$ are invariants of minimum weight.
It follows from Lemma~\ref{45dellem} that  $\delta^{2p-d-3}(x_4^{d-p+1}x_5^{p-1})$ is a linear combination
of these invariants with the coefficient of $x_4^{p-2}x_5^{d-p+2}$ non-zero.
It follows from Lemma~\ref{345dellem} that $\delta^{2p-d-4}(x_3x_4^{d-p}x_5^{p-1})$ is a linear
combination of $x_3x_4^{p-4}x_5^{d-p+3}$ and $x_2x_4^{p-5}x_5^{d-p+4}$
with the coefficient of $x_3x_4^{p-4}x_5^{d-p+3}$ non-zero.
Since $x_2$ is invariant, $\delta^{2p-d-5}(x_2x_4^{d-p}x_5^{p-1})=x_2\delta^{2p-d-5}(x_4^{d-p}x_5^{p-1})$
which, by Lemma~\ref{45dellem}, is a non-zero scalar multiple of $x_2x_4^{p-5}x_5^{d-p+4}$.
Thus we have submodules of dimensions $2p-d-2$, $2p-d-3$ and $2p-d-4$ such that the sum of the socles is direct.
By Lemma~\ref{soclelem}, this gives the required
decomposition.

(iv)($d=p,\,p+1$) As in (iii), the elements $x_4^{d-p+1}x_5^{p-1}$, $x_3x_4^{d-p}x_5^{p-1}$ and $x_2x_4^{d-p}x_5^{p-1}$
generate submodules of dimensions  $2p-d-2$, $2p-d-3$ and $2p-d-4$, respectively.
Furthermore, $\delta^{2p-d-3}(x_4^{d-p+1}x_5^{p-1})$, $\delta^{2p-d-4}(x_3x_4^{d-p}x_5^{p-1})$ and
$\delta^{2p-d-5}(x_2x_4^{d-p}x_5^,{p-1})$ are linearly independent elements of 
${\rm Span}(x_4^{p-2}x_5^{d-p+2},x_3x_4^{p-4}x_5^{d-p+3}, x_2x_4^{p-5}x_5^{d-p+4})$.
The invariant basis element  $x_2x_3x_5^{d-2}$ generates a submodule of dimension $1$.
Thus the sum of the socles of these four submodules is direct and the sum of their dimensions
is the dimension of the homogeneous component. Therefore, by Lemma~\ref{soclelem}, we have the required decomposition.

(v) ($d=p-2,\,p-1$) As in (iii), the  elements  $x_4^{p-2}x_5^{d-p+2}$,
$x_3x_4^{p-4}x_5^{d-p+3}$ and $x_2x_4^{p-5}x_5^{d-p+4}$ are invariants of minimum weight.
Using Lemma~\ref{45dellem} and Lemma~\ref{345dellem}, the basis elements
$x_5^d$, $x_3x_5^{d-1}$ and $x_2x_5^{d-1}$ generate submodules of dimensions
$p-1$, $p-3$ and $p-4$, respectively. Furthermore, $\delta^{p-1}(x_5^d)$,  $\delta^{p-3}(x_3x_5^{d-1)}$
and $\delta^{p-4}(x_2x_5^{d-1})$ are linearly independent minimum weight invariants.
As in (iv), the invariant basis element  $x_2x_3x_5^{d-2}$ is a non-minimum weight invariant.
Thus, applying Lemma~\ref{soclelem}, we have four submodules whose sum is direct and whose dimensions sum to the dimension of the homogeneous component.

(vi) ($d=p-3$, $p>11$) From Lemma~\ref{numlem}, $\vcoin^{p-3}$ is a sum of at most four indecomposable summands.
We will identify four linearly independent invariants. Since each summand has a one dimensional
socle, this means that there are four summands  and we have found a basis for the invariants.
In this homogeneous component, the minimum weight is $4p-13$ and the minimum weight subspace is 
${\rm Span}(x_3x_4^{p-4},x_2x_4^{p-5})$. This gives two linearly independent invariants.
 The weight $4p-12$ subspace is ${\rm Span}(x_4^{p-3},x_3x_4^{p-5},x_2x_4^{p-6}x_5^2)$.
(Note that $\wt(x_2x_3x_5^{p-5})=5p-20$. Therefore, since $p>8$, $\wt(x_2x_3x_5^{p-5})>4p-12$.)
Since $\delta$ is a linear map taking the three dimensional weight $4p-12$ subspace to the two dimensional
weight $4p-13$ subspace, there exists a non-zero element 
$f\in {\rm Span}(x_4^{p-3},x_3x_4^{p-5},x_2x_4^{p-6}x_5^2)$ with $\delta(f)=0$. Using Lemma~\ref{wtlem},
$\Delta(f)=\delta(f)$. Therefore $f$ is invariant. The fourth invariant is $x_2x_3x_5^{p-3}$.

From Lemma~\ref{5dellem}(i), 
$$\delta^{p-2}(x_5^{p-3})=\frac{(p-3)(p-2)!}{12}(6x_3x_4^{p-4}+(p-4)(3p-13)x_2x_4^{p-5}x_5),$$
while from Lemma~\ref{345dellem},
$$\delta^{p-4}(x_3 x_5^{p-4})=(p-4)!\left(x_3 x_4^{p-4} +(p-4)^2 x_2x_4^{p-5}x_5\right).$$
A simple calculation shows that for $p>11$,  $6(p-4)^2\not\equiv_{(p)}(p-4)(3p-13)$.
Therefore $\delta^{p-4}(x_3 x_5^{p-4})$  and $\delta^{p-2}(x_5^{p-3})$ are linearly independent
minimum weight invariants. Using Lemma~\ref{45dellem}, $\delta^{p-5}(x_2x_5^{p-4})=(p-4)!x_2x_4^{p-5}x_5$.
Thus $x_2x_5^{p-4}$ generates a submodule of dimension $p-4$ whose socle is contained in the
the minimum weight subspace. Since  $\delta^{p-4}(x_3 x_5^{p-4})$  and $\delta^{p-2}(x_5^{p-3})$
are a basis for the minimum weight subspace, it is possible to choose coefficients $c_1$ and $c_2$
so that $\delta^{p-5}(x_2x_5^{p-4}+c_1\delta(x_3 x_5^{p-4})+c_2 \delta^3(x_5^{p-3}))=0$.
We claim that $h:=\delta^{p-4}(x_2x_5^{p-4}+c_1\delta(x_3 x_5^{p-4})+c_2 \delta^3(x_5^{p-3}))$ is
 non-zero scalar multiple of $f$. Clearly $h$ is an invariant of weight $4p-12$.
However, if $h$ is zero, then $f$ is not contained in the submodule generated by
$\{x_5^{p-3},x_3x_5^{p-4},x_2x_5^{p-4},x_2x_3x_5^{p-2}\}$ contradicting Lemma~\ref{numlem}.

In conclusion, $x_5^{p-3}$, $x_3x_5^{p-4}$, $x_2 x_5^{p-4}+c_1\delta(x_3 x_5^{p-4})+c_2 \delta^3(x_5^{p-3})$,
and $x_2x_3x_5^{p-4}$ generate submodules of dimensions $p-1$, $p-3$, $p-5$ and $1$, respectively.
The sum of the socles of these modules is direct and the sum of the dimensions matches the dimension of
$\vcoin^{p-3}$. Therefore, by Lemma~\ref{soclelem}, $\vcoin^{p-3}\cong V_{p-1}\oplus V_{p-3} \oplus V_{p-5}\oplus V_1$.

(vii) For $d=6,\ldots, p-4$ a basis for $\vcoin^d$ is given by
$$
\begin{array}{llllllll}
               &                  &x_4^d,              &\ldots\ldots& x_4^3 x_5^{d-3},  &x_4^2 x_5^{d-2}, & x_4 x_5^{d-1}, & x_5^d,\\
               &x_3 x_4^{d-1},    &x_3 x_4^{d-2} x_5,  &\ldots\ldots& x_3 x_4 x_5^{d-2},&x_3 x_5^{d-1},\\
x_2 x_4^{d-1}, &x_2 x_4^{d-2}x_5, &x_2 x_4^{d-3} x_5^2, &\ldots\ldots& x_2 x_5^{d-1},\\ 
&&&x_2 x_3 x_5^{d-2},
\end{array}
$$
where elements in the same column have the equal weight.
The case $d=5$ is essentially the same except the element in the fourth row lies in the third column.
As in (vi), using Lemma~\ref{numlem}, $\vcoin^d$ is the sum of at most four indecomposable summands.
We first show that there are four linearly independent invariants and hence four summands.
The elements $x_2 x_4^{d-1}$, $x_2 x_3 x_5^{d-2}$, and $x_3 x_4^{d-1} - (2d-1)x_2 x_4^{d-2}x_5$
are easily seen to be invariant and linearly independent. A fourth invariant
can be constructed as a linear combination of $x_4^d$, $x_3 x_4^{d-2} x_5$, $x_2 x_4^{d-3} x_5^2$
and $x_2 x_4^{d-2}x_5$. To see this, first observe that the weight $4d$ subspace has dimension 3 ($4$ for $d=5$)
and that the weight $4d-1$ subspace has dimension 2. Therefore there is a non-zero linear combination
of $x_4^d$, $x_3 x_4^{d-2} x_5$ and $x_2 x_4^{d-3} x_5^2$ in the kernel of $\delta$, say $f$. By Lemma~\ref{wtlem}, $\delta(f)$ and 
$\del(f)$ agree in weight $4d-1$. Furthermore, the only basis element of lower weight is $x_2 x_4^{d-1}$. Thus 
$\del(f)$ is a scalar multiple of $x_2 x_4^{d-1}$. Note that $\del(x_2 x_4^{d-2}x_5)$ is a non-zero scalar multiple
of  $x_2 x_4^{d-1}$. Therefore there exists $a\in\field$ with $\del(f-a x_2 x_4^{d-2}x_5)=0$, giving
the required invariant.
Thus we have four linearly independent invariants and four indecomposable summands.

Suppose $3d-1\not\equiv_{(p)} 0$.  From Lemma~\ref{5dellem}(ii),
$\delta^{d+2}(x_5^d)=x_3 x_4^{d-1}d(3d-1)(d+2)!/12$. Thus $x_5^d$ generates a submodule of dimension $d+3$ 
with socle ${\rm Span}(x_2 x_4^{d-1})$. From 
Lemma~\ref{35dellem}, $\delta^d(x_3 x_5^{d-1})=d(d!)x_2 x_4^{d-1}$. Thus 
$\delta^d(x_3 x_5^{d-1}-c\delta^2(x_5^d))=0$ with $c=12/((d+1)(d+2)(3d-1))$.
Using Lemma~\ref{345dellem} and Lemma~\ref{5dellem}(i), we have
$$\begin{array}{lll}
\delta^{d-1}\left(x_3 x_5^{d-1}-c\delta^2(x_5^d)\right)=&\left((d-1)!-\frac{cd(d+1)!}{2}\right)x_3 x_4^{d-1}\\
&\\
& +\left((d-1)^2(d-1)!-\frac{cd(d-1)(3d-4)(d+1)!}{12}\right)x_2 x_4^{d-2} x_5.
  \end{array}$$
Substituting for $c$ and simplifying gives
$$\delta^{d-1}\left(x_3 x_5^{d-1}-c\delta^2(x_5^d)\right)=\frac{-(d-1)(3d-2)(d-1)!}{(d+2)(3d-1)}
\left(x_3 x_4^{d-1} + (2d-1)x_2 x_4^{d-2}x_5\right).$$
Note that $\delta^{d-1}\left(x_3 x_5^{d-1}-c\delta^2(x_5^d)\right)$ and
$\del^{d-1}\left(x_3 x_5^{d-1}-c\delta^2(x_5^d)\right)$ differ by a scalar multiple of
the invariant $x_2x_4^{d-1}$.
Suppose $3p-2\not\equiv_{(p)}0$. Then $x_3 x_5^{d-1}-c\delta^2(x_5^d)$ generates a module of 
dimension $d$. Furthermore, it is possible to choose a linear combination of
$x_2 x_5^{d-1}$, $x_3 x_4 x_5^{d-2}$, and $x_4^3 x_5^{d-3}$, say $h$, so that $h$ generates a submodule
of dimension $d-2$ with socle given by the span of $f+a x_2 x_4^{d-2}x_5$. Thus applying
Lemma~\ref{soclelem} gives a decomposition isomorphic to $V_{d+3}\oplus V_d \oplus V_{d-2}\oplus V_1$.
On the other hand, suppose $3p-2\equiv_{(p)}0$. Then $\del^{d-1}\left(x_3 x_5^{d-1}-c\delta^2(x_5^d)\right)$
is a multiple of $x_2 x_4^{d-1}$ and, for some $c'\in\field$, $x_3 x_5^{d-1}-c\delta^2(x_5^d)+c'x_4^3 x_5^{d-3}$
generates a module of dimension at most $d-1$. Furthermore, it is possible to choose a linear combination of
$x_2 x_5^{d-1}$, $x_3 x_4 x_5^{d-2}$, and $x_4^3 x_5^{d-3}$, say $h'$, so that
$h'$ generates a submodule of dimension $d-1$ with socle given by the span of $x_3 x_4^{d-1} + (2d-1)x_2 x_4^{d-2}x_5$.
Since $\{x_5^d, x_3 x_5^{d-1}, x_2 x_5^{d-1}, x_2 x_3 x_5^{d-2}\}$ is generating set for the homogeneous component,
the module generated by $x_3 x_5^{d-1}-c\delta^2(x_5^d)+c'x_4^3 x_5^{d-3}$ has dimension $d-1$ and its
socle does not lie in ${\rm Span}\left(x_2 x_4^{d-1}, x_2 x_3 x_5^{d-2},x_3 x_4^{d-1} + (2d-1)x_2 x_4^{d-2}x_5\right)$.
Thus, using Lemma~\ref{soclelem}, we have a decomposition isomorphic to
$V_{d+3}\oplus 2V_{d-1}\oplus V_1$.

Suppose $3d-1\equiv_{(p)}0$. By Lemma~\ref{5dellem}, $x_5^d$ generates a module of dimension $d+2$
with socle ${\rm Span}(3 x_3 x_4^{d-1}+ x_2 x_4^{d-2} x_5+cx_2 x_4^{d-1})$ for some $c\in\field$. By Lemma~\ref{35dellem}
$x_3 x_5^{d-1}$ generates a module of dimension $d+1$ with socle ${\rm Span}(x_2 x_4^{d-1})$.
Clearly $x_2 x_3 x_5^{d-2}$ generates a module of dimension $1$.
Since $\{x_5^d, x_3 x_5^{d-1}, x_2 x_5^{d-1}, x_2x_3 x_5^{d-2} \}$ generates the homogeneous component,
a suitable linear combination of
$x_2 x_5^{d-1}$, $x_3 x_4 x_5^{d-2}$, and $x_4^3 x_5^{d-3}$ generates a module of
dimension $d-2$ with socle determined by $f-a x_2 x_4^{d-2}x_5$. Applying Lemma~\ref{soclelem} gives
a decomposition isomorphic to $V_{d+2}\oplus V_{d+1}\oplus V_{d-2}\oplus V_1$.

(viii) For $d=3$: The dimension of the homogeneous component is $11$. A direct calculation of
$\delta^5(x_5^3)$ shows that $x_5^3$ generates a module of dimension $6$ with
with socle ${\rm Span}(4x_2 x_4^2 +x_2 x_3 x_5)$. A direct calculation of $\delta^3(x_3 x_5^2)$
shows that $x_3 x_5^2$ generates a submodule of dimension $4$ with socle ${\rm Span}(9x_2 x_4^2+4 x_2 x_3 x_5)$.
The linear map $\del$ takes the span of the elements of weight less than $13$, a subspace of dimension $7$, to
the the span of the elements of weight less than $12$, a subspace of dimension $4$. Thus the kernel of
$\del$ has dimension at least $3$. Applying Lemma~\ref{soclelem} gives the required decomposition

For $d=4$: It is clear that $x_2 x_4^3$, $x_3 x_4^3 - 7 x_2 x_4^2 x_5$ and $x_2 x_3 x_5^2$ are invariant.
Note that, for $p>7$, the dimension of the homogeneous component is $14$.
From Lemma~\ref{5dellem}, $\delta^6(x_5^4)=2640 x_2 x_4^3$. Thus, for $p>11$, $x_5^4$ generates a module of dimension
$7$ with socle ${\rm Span}(x_2 x_4^3)$. From Lemma~\ref{35dellem}, $\delta^4(x_3 x_5^3)= 96 x_2 x_4^3$.
Define $g_1:=55 x_3 x_5^3-2\delta^2(x_5^4)$. Using Lemma~\ref{5dellem}(i)  and Lemma~\ref{345dellem}, 
$\delta^3(g_1)=-150(x_3 x_4^3 - 7 x_2 x_4^2 x_5)+420x_2 x_3 x_5^2$. Thus $g_1$ generates a module of
dimension $4$. Using Lemma~\ref{45dellem}, $\delta^3(x_2 x_5^3)= 6 x_2 x_4^3$. Define
$g_2:=440 x_2 x_5^3-\delta^3(x_5^4)$. From Lemma~\ref{45dellem}, 
$\delta^2(g_2)=-240( x_3 x_4^3 - 7 x_2 x_4^2 x_5)+1200 x_2 x_3 x_5^2$.
Thus $g_2$ generates a module of dimension $3$.
For $p>11$, $\delta^3(g_2)$ and $\delta^4(g_1)$ are linearly independent. Thus  $\del^3(g_2)$ and $\del^4(g_1)$
are linearly independent and applying Lemma~\ref{soclelem} gives the required decomposition.
\end{proof}

\end{document}